%% file: hermite.tex
\documentclass[11pt]{article}
\input{macros_fr_pdf.tex}

\title{Sur les minima des formes hamiltoniennes binaires \\
  d\'efinies positives}
\author{Ga\"etan Chenevier \and Fr\'ed\'eric Paulin} 
\date{\today}

\begin{document}
\bibliographystyle{alphanum}
\maketitle
\begin{abstract}\footnote{{\bf Mots clefs :} alg\`ebre de quaternions, forme
    hamiltonienne binaire, ordre maximal, r\'eseau euclidien.~~
    {\bf Codes AMS :} 11E39, 11R52, 11L05, 16H20, 11E20.}
\'Etant donn\'e un ordre maximal $\OOO$ d'une alg\`ebre de quaternions
rationnelle d\'efinie $A$ de discriminant $D_A$, nous montrons que le
minimum des formes hamiltoniennes binaires sur $\OOO$, d\'efinies
positives et de discriminant $-1$, est $\sqrt{D_A}$.  Lorsque la
diff\'erente de $\OOO$ est principale, nous explicitons une forme
atteignant cette valeur, et lorsque $\OOO$ est principal, nous donnons
la liste exacte des formes atteignant cette valeur.  Nous donnons des
crit\`eres et des algorithmes pour d\'eterminer quand la diff\'erente de
$\OOO$ est principale. 

\medskip
Let $A$ be a definite quaternion algebra over $\QQ$, with discriminant $D_A$, and $\OOO$ a maximal order of $A$.
We show that the minimum of the positive definite hamiltonian binary forms over $\OOO$ with discrimiminant $-1$ is $\sqrt{D_A}$.
When the different of $\OOO$ is principal, we provide an explicit form representing this minimum, and when $\OOO$ is principal, 
we give the list of the equivalence classes of all such forms. 
We also give criteria and algorithms to determine when the different of $\OOO$ is principal.
\end{abstract}

\section{Introduction}
\label{subsect:introB}

Soit $A$ une alg\`ebre de quaternions sur $\QQ$ qui est d\'efinie, de
sorte que $\HH=A\otimes_\QQ\RR$ soit l'alg\`ebre des quaternions de
Hamilton sur $\RR$ usuelle. Nous noterons $x \mapsto \overline{x}$ la
conjugaison, $\n$ la norme r\'eduite, $\tr$ la trace r\'eduite de
$\HH$, $D_A$ le discriminant r\'eduit de $A$, $h_A$ son nombre de
classes et $t_A$ son nombre de classes de conjugaison d'ordres
maximaux. Soit $\OOO$ un ordre maximal dans $A$. Rappelons que la
{\it diff\'erente} de $\OOO$ est l'unique id\'eal \`a droite 
de $\OOO$ de norme (r\'eduite) $D_A$, et qu'il est bilat\`ere. Nous renvoyons par exemple
\`a \cite{Vigneras80} pour les pr\'erequis.

En utilisant la terminologie de \cite{Weyl40}, notons $\Q^+$
l'ensemble des {\it formes hamiltoniennes binaires} 
$$
f :(u,v)\mapsto a\n(u)+ \tr(\overline{u}\,b\,v) +c\n(v)
$$ 
(o\`u $b,u,v\in \HH$ et $a,c\in\RR$) qui sont d\'efinies positives, ou
de mani\`ere \'equivalente avec $a,c>0$ et de {\it discriminant}
$\Delta(f)=\n(b)-ac$ strictement n\'egatif, et $\Q^+_1$ le
sous-ensemble des $f\in \Q^+$ telles que $\Delta(f)=-1$.
D\'efinissons la {\it constante d'Hermite $\ga_2(\OOO)$ de l'ordre
  maximal $\OOO$} par
$$
\ga_2(\OOO)=\sup_{f\in\Q^+}\;\min_{(u,v)\in\OOO\times\OOO-\{0\}}
\;\frac{1}{\sqrt{-\Delta(f)}}\;f(u,v)=
\sup_{f\in\Q^+_1}\;\min_{(u,v)\in\OOO\times\OOO-\{0\}}\;f(u,v)\;.  
$$

Si nous rempla\c{c}ons $\OOO$ par $\ZZ$ et si $f$ varie parmi les formes
quadratiques binaires r\'eelles d\'efinies positives, cette constante
$\ga_2=\ga_2(\ZZ)$, appel\'ee la constante d'Hermite binaire, vaut
$\frac{2}{\sqrt{3}}$, et la description des formes $f$ qui r\'ealisent
la borne sup\'erieure est bien connue (voir par exemple
\cite[p.~332]{Cassels59}).  Nous renvoyons \`a \cite{Oppenheim34} pour 
le cas o\`u $\OOO$ est remplac\'e par l'anneau des entiers d'une extension
quadratique imaginaire de $\QQ$, quand $f$ varie sur les formes
hermitiennes binaires complexes d\'efinies positives
(par exemple $\ga_2(\ZZ[i]) =2$).

Le r\'esultat principal de cette note, g\'en\'eralisant le cas
$D_A=2$ trait\'e par \cite[Satz 4]{Speiser32}, est le suivant.

\btheo\label{theo:appBminbinhamform}
Nous avons
$
\ga_2(\OOO)= \sqrt{D_A}\;.
$
\etheo

L'in\'egalit\'e $\ga_2(\OOO)\leq \sqrt{D_A}$ d\'ecoulera assez
facilement du calcul (voir \cite{Blichfeldt35}) de la constante de
Hermite $\ga_8$ pour les formes quadratiques r\'eelles d\'efinies
positives en $8$ variables.  L'\'egalit\'e r\'esultera du fait que le
r\'eseau euclidien ${\rm E}_8$ peut \^etre muni, pour tout ordre
maximal $\OOO$, d'une structure de $\OOO$-module libre 
(et pas seulement projectif) de rang $2$, pour
laquelle les \'el\'ements $x$ de $\OOO$ agissent par des similitudes
orthogonales de rapport $\n(x)$ (voir la proposition \ref{prop:E8}).
Notre construction, reposant sur des techniques de r\'esidus de r\'eseaux euclidiens, 
g\'en\'eralise une construction classique de ${\rm E}_8$
\`a l'aide des quaternions de Hurwitz (voir par exemple \cite[Prop.~8.2.2]{Martinet03}).
Cette g\'en\'eralisation est tr\`es directe lorsque la diff\'erente de $\OOO$
est suppos\'ee principale, et conduit dans ce cas \`a des 
descriptions explicites : voir la proposition \ref{prop:Lambdalambda}. 
\medskip

L'in\'egalit\'e $\ga_2(\OOO)\leq \sqrt{D_A}$ est utilis\'ee dans
\cite{ParPau19}, qui donne \`a l'aide d'outils de g\'eom\'etrie
hyperbolique de dimension $5$ une th\'eorie graphique des formes
hamiltoniennes binaires enti\`eres ind\'efinies, analogue \`a celle de
Conway pour les formes quadratiques binaires.

\medskip
Dans la partie \ref{sec:casegal}, lorsque la diff\'erente de $\OOO$ est
principale, nous donnons une liste (a priori incompl\`ete, certainement redondante) de formes hamiltoniennes binaires d\'efinies positives
atteignant $\ga_2(\OOO)$. Lorsque $\OOO$ lui-m\^eme est principal
(c'est-\`a-dire lorsque $h_A=1$), nous \'etudions l'unicit\'e d'une telle
forme. Deux formes hamiltoniennes binaires seront
dites {\it $\OOO$-\'equivalentes} si elles se d\'eduisent l'une de 
l'autre par pr\'ecomposition par un \'el\'ement de $\GL_2(\OOO)$.

\bprop \label{prop:nonuniintro} Si $D_A =2,3,5,7$, il existe une
unique classe de $\OOO$-\'equivalence 
de formes hamiltoniennes binaires d\'efinies positives de discriminant donn\'e et r\'ealisant la
borne sup\'erieure d\'efinissant la constante d'Hermite
$\ga_2(\OOO)$. Si $D_A=13$, il en existe exactement deux.  
\eprop

Une \'etude g\'en\'erale des classes de $\OOO$-\'equivalence de telles
formes serait int\'eressante (voir la remarque \ref{remxo}).
Nous revenons dans la partie \ref{subsect:ordreetreseaux} sur la
correspondance classique entre classes de conjugaison d'ordres maximaux de $A$
et certaines formes quadratiques ternaires, par des techniques de r\'esidus et
sommes de Gauss. Cette correspondance associe \`a l'ordre maximal $\OOO$
deux r\'eseaux euclidiens pairs de dimension $3$ :
d'une part  $L(\OOO)=\{x\in
\OOO\;:\; \tr x=0\}$ muni de la restriction de la forme norme, et d'autre part le plus grand sous-r\'eseau pair $M(\OOO)$ de $N \cap N^\sharp$ avec $N =\frac{1}{\sqrt{D_A}} L(\OOO)$.
Dans la partie \ref{subsect:ordmaxdiffprinc} (voir le th\'eor\`eme \ref{theo:carac1}), nous montrons alors les \'equivalences entre~: 
\medskip

$\bullet$ la diff\'erente de $\OOO$ est principale, 

$\bullet$ $\OOO$ contient un \'el\'ement de carr\'e $-D_A$, 

$\bullet$ $L(\OOO)$ contient un \'el\'ement $x$ tel
que $x\cdot x = 2D_A$ et $x\cdot y = 0 \mod D_A$ pour tout $y \in L(\OOO)$,
et

$\bullet$~ $M(\OOO)$ contient un \'el\'ement $x$ tel que $x \cdot x = 2$.

\medskip
\noindent 
Nous utilisons ces \'equivalences pour donner
de nombreux exemples (voir la proposition \ref{prop:tdnpstricpositf}
montrant que $A$ admet toujours au moins un ordre maximal de
diff\'erente principale) et contre-exemples (voir le tableau final de
cette note).

\medskip\noindent{\small {\it Remerciements : } Le premier auteur a
  \'et\'e financ\'e par le C.N.R.S.  et a re\c{c}u le soutien du
  projet ANR-14-CE25 (PerCoLaTor). Le second auteur remercie
  l'universit\'e de Warwick et l'EPSRC pour leur accueil et soutien
  financier lors de la r\'edaction d'une premi\`ere version de cette
  note. Les auteurs remercient Lassina Demb\'el\'e d'avoir v\'erifi\'e certains de leurs calculs avec le logiciel Magma.  }

\section{Calcul de la constante d'Hermite binaire $\ga_2(\OOO)$}

\label{subsect:demomajo}

Dans toute cette note, nous munissons $\HH$ du produit scalaire
euclidien $(x,y)\mapsto \tr (\,\overline{x}\,y)$ (rendant la base
usuelle $(1,i,j,k)$ de $\HH$ orthogonale et constitu\'ee de vecteurs
de norme $\sqrt{2}$), et $\HH\times\HH$ du produit scalaire euclidien
produit, que nous notons $(w,w')\mapsto w\cdot w'$. En particulier,
ceci d\'efinit une forme volume sur $\HH\times\HH$, et nous
d\'efinissons le {\it covolume} d'un $\ZZ$-r\'eseau $\Lambda$ de
$\HH\times\HH$ par
$$ 
\covol\Lambda=\Vol((\HH\times\HH)/\Lambda)\;.  
$$
Puisque $\OOO$ est un $\ZZ$-r\'eseau dans $\HH$, le produit
$\OOO\times\OOO$ est un $\ZZ$-r\'eseau dans $\HH\times\HH$, de covolume
(voir \cite[Lem.~5.5]{KraOse90})
$$
\covol(\OOO\times\OOO)=\Vol(\HH/\OOO)^2= D_A^{\;\;\,2}\;.  
$$

Nous munissons $\HH\times\HH$ de sa structure
d'espace vectoriel \`a droite sur $\HH$. Soit $f_0$ la forme hamiltonienne binaire d\'efinie positive
$$
f_0(u,v)=\n(u)+\n(v)\;,
$$ 
de discriminant $-1$, de sorte que pour tout $w\in\HH\times\HH$, nous
avons $f_0(w)=\frac{1}{2}\,w\cdot w$. Nous noterons $h$ le produit scalaire hermitien sur $\HH \times \HH$ 
tel que $h(w,w)=f_0(w)$ pour tout $w \in \HH \times \HH$. Pour tout $w'=(x',y') \in \HH \times \HH$ et tout $w = (x,y) \in \HH \times \HH$, 
nous avons $h(w',w)=\overline{x}x'+\overline{y}y'$.

L'action \`a droite par pr\'ecomposition du groupe $\SL_2(\HH)$ des
matrices $2\times 2$ \`a coefficients dans $\HH$ de d\'eterminant de
Dieudonn\'e $1$ est transitive sur $\Q^+_1$ (voir par exemple \cite[\S
  7]{ParPau13ANT}). On appelle
  {\em $\OOO$-r\'eseau} de $\HH \times \HH$ un $\ZZ$-r\'eseau qui est en outre un sous-module {\it libre} de rang $2$
  du $\OOO$-module \`a droite $\HH \times \HH$. L'action lin\'eaire \`a
gauche de $\SL_2(\HH)$ sur l'ensemble des $\OOO$-r\'eseaux de $\HH\times\HH$ qui
sont de covolume donn\'e en tant que $\ZZ$-r\'eseaux est aussi
transitive.

Pour tout $n\in\NN-\{0\}$, rappelons (voir par exemple
\cite{Cassels59}) que, avec $Q^+(n)$ l'ensemble des formes
quadratiques r\'eelles d\'efinies positives en $n$ variables, $M(f)$
la matrice\footnote{de sorte que $f(x)=\;^t\!XM(f)X$ si $X$ est la
  matrice colonne des coordonn\'ees de $x$} de $f$ et $(x,y)\mapsto
x\cdot y$ le produit scalaire usuel sur $\RR^n$, la {\it constante
  d'Hermite} $\ga_n$ en dimension $n$ est d\'efinie par
$$
\ga_n=\sup_{f\in Q^+(n)}\;
\min_{x\in\ZZ^n-\{0\}}\;\frac{f(x)}{\sqrt[n]{\det M(f)}}
\;\;=\;\sup_{\substack{L \;\;\ZZ\text{-r\'eseau~de}\;\RR^n\\ \covol L\,=\,1}}\;
\min_{x\in L-\{0\}}\; x\cdot x\;.
$$
Les valeurs de  $\ga_n$ sont connues si $n\leq 8$, par exemple  Blichfeld
\cite{Blichfeldt35} a montr\'e que
$$
\ga_8=2\;.
$$
Il est bien connu que cette valeur est atteinte pour le r\'eseau ${\rm
  E}_8$, engendr\'e par un syst\`eme de racines de longueur $2$ et de type $E_8$, 
 qui contient $240$ vecteurs $x$ tels que $x\cdot x=2$. On rappelle que d'apr\`es
 Mordell \cite{Mordell38}, 
 ${\rm E}_8$ est l'unique (\`a isom\'etrie pr\`es) $\ZZ$-r\'eseau euclidien
entier pair, unimodulaire ({\it i.e.} de covolume $1$) et de rang $8$ : voir les rappels
ci-dessous pour ces terminologies classiques. Mieux, on 
sait d'apr\`es Vetchinkin \cite[Theo.~2]{Vetchinkin80} qu'\`a
isom\'etrie pr\`es, ${\rm E}_8$ est le seul $\ZZ$-r\'eseau unimodulaire
atteignant la borne sup\'erieure d\'efinissant $\ga_8$.

Tout $\OOO$-r\'eseau de $\HH\times
\HH$ est en particulier un $\ZZ$-r\'eseau de l'espace euclidien r\'eel $\HH\times\HH$
de dimension $8$. Donc
\begin{align*}
\ga_2(\OOO)&=\sup_{f\in\Q^+_1}\;\min_{w\in\OOO\times\OOO-\{0\}}\;f(w)
=\sup_{g\in\SL_2(\HH)}\;\min_{w\in \OOO\times\OOO-\{0\}}\;f_0\circ g(w)\\ &
=\sup_{g\in\SL_2(\HH)}\;\min_{w\in g(\OOO\times\OOO)-\{0\}}\;f_0(w)
=\sup_{\substack{\Lambda \;\;\OOO\text{-r\'eseau~de}\;\HH\times\HH\\ 
\covol \Lambda\,=\,D_A^2}}\;\;
\min_{w\in \Lambda-\{0\}}\; w\cdot w/2\\ &
=\sqrt{D_A}\;
\sup_{\substack{\Lambda \;\;\OOO\text{-r\'eseau~de}\;\HH\times\HH\\ 
\covol \Lambda\,=\,1}}\;\;
\min_{w\in \Lambda-\{0\}}\; w\cdot w/2\\ &
\leq\frac{\sqrt{D_A}}{2}\;
\sup_{\substack{L \;\;\ZZ\text{-r\'eseau~de}\;\RR^8\\ \covol L\,=\,1}}\;
\min_{x\in L-\{0\}}\; x\cdot x \; = \frac{\sqrt{D_A} }{2}\;\ga_8=\sqrt{D_A}
\;. 
\end{align*}
De plus, ce calcul montre que nous avons \'egalit\'e dans
l'in\'egalit\'e $\ga_2(\OOO)\leq \sqrt{D_A}$ si et seulement s'il
existe un $\OOO$-r\'eseau $\Lambda$ de covolume $1$ dans
$\HH\times\HH$ tel que $\min_{w\in \Lambda-\{0\}}\; w\cdot w=2$.  Le
th\'eor\`eme \ref{theo:appBminbinhamform} d\'ecoule donc du r\'esultat
suivant.

\bprop \label{prop:E8} L'espace euclidien $\HH \times \HH$ contient
des $\OOO$-r\'eseaux isom\'etriques \`a ${\rm E}_8$.
\eprop

Nous allons utiliser la technique bien connue des r\'esidus de
$\ZZ$-r\'eseaux euclidiens (aussi appel\'es ``formes quadratiques
discriminantes'', voire ``glue groups'' dans \cite{ConSlo88}), voir
par exemple \cite[\S 3.3]{Ebeling13}, \cite[\S II.1]{CheLan19}.
Rappelons bri\`evement les \'el\'ements utiles \`a notre propos.

Un {\it module quadratique d'enlacement}, ou pour faire court un ${\rm
  qe}$-{\it module} au sens de \cite[\S II.1]{CheLan19}, est la
donn\'ee d'un groupe ab\'elien fini $V$ muni d'une application $q: V
\rightarrow \QQ/\ZZ$ v\'erifiant $q(mx)=m^2q(x)$ pour tous les
$m\in\ZZ$ et $x$ dans $V$, et telle que l'application $V \times V
\rightarrow \QQ/\ZZ$, d\'efinie par $(x,y)\mapsto q(x+y)-q(x)-q(y)$, est
$\ZZ$-bilin\'eaire non d\'eg\'en\'er\'ee. Un ${\rm qe}$-{\it module}
$(V,q)$ est {\it anisotrope} si le seul \'el\'ement $x$ de $V$ tel que
$q(x)=0$ est $x=0$.

Soient $E$ un espace vectoriel r\'eel euclidien de produit scalaire
not\'e $(x,y)\mapsto x\cdot y$, et $L$ un $\ZZ$-r\'eseau de $E$.
Notons $L^\sharp=\{x\in E\;:\; \forall\;y\in L,\;\;x\cdot y\in\ZZ\}$
le $\ZZ$-r\'eseau {\it dual} de $L$. Le r\'eseau $L$ est dit {\it entier} si l'on a 
$L \subset L^\sharp$, {\it pair} si de plus $x \cdot x \in 2\ZZ$ pour tout $x \in L$. 
Supposons d\'esormais $L$ entier et
pair. Alors l'application
\begin{equation}\label{eq:res}
  q: L^\sharp/L \longrightarrow \QQ/\ZZ,
  \hspace{5mm} x+L\mapsto \frac{x\cdot x}{2}\mod \ZZ,
\end{equation}
munit $L^\sharp/L$ d'une structure de ${\rm qe}$-module not\'ee $\res
L$ et appel\'ee le {\it r\'esidu} de $L$.

Le groupe ab\'elien fini $L^\sharp/L$ est d'ordre le {\it
  d\'eterminant} de $L$ (le d\'eterminant de la matrice de Gram de
n'importe quelle $\ZZ$-base de $L$), not\'e $\det L$ ; nous avons
$$
\det L' = [L:L']^2 \det L
$$ 
si $L'$ est un sous-$\ZZ$-r\'eseau de $L$. La projection canonique
$\operatorname{pr} : L^\sharp\ra L^\sharp/L$ induit une bijection de
l'ensemble des r\'eseaux entiers et pairs $\Lambda$ contenant $L$ avec
indice $m$ sur l'ensemble des sous-groupes $I$ de $\res L$, d'ordres
$m$ et {\em isotropes} (tels que $q(x)=0$ pour tout $x\in I$), de
sorte que $\Lambda=\operatorname{pr}^{-1}(I)$.

Pour tout nombre premier $p$, on pose $L_p = L \otimes_\ZZ \ZZ_p$.
C'est un $\ZZ_p$-module muni de la forme $\ZZ_p$-bilin\'eaire $L_p
\times L_p \rightarrow \ZZ_p, \, (x,y) \mapsto x \cdot y,$ d\'eduite
par extension des scalaires du produit scalaire sur $L$.  On a $x
\cdot x \in 2\ZZ_p$ pour tout $x$ dans $L_p$. Le $\ZZ_p$-r\'eseau
$L_p$ poss\`ede un {\em dual} d\'efini par $L_p^\sharp = \{x\in
L_p[\frac{1}{p}] \;:\; \forall\;y\in L_p,\;\;x\cdot y\in\ZZ_p\}$,
contenant $L_p$, ainsi qu'un {\em r\'esidu} $\res L_p$ qui est
le $p$-groupe ab\'elien fini $L_p^\sharp/L_p$ muni de la forme
quadratique \`a valeurs dans $\QQ_p/\ZZ_p$ d\'efinie comme dans
\eqref{eq:res}.  L'application naturelle $\QQ/\ZZ \rightarrow
\QQ_p/\ZZ_p$ induit un isomorphisme de la composante $p$-primaire du
groupe ab\'elien de torsion $\QQ/\ZZ$ vers $\QQ_p/\ZZ_p$.  Cet
isomorphisme permet de voir $\res L_p$ comme un ${\rm
  qe}$-module.  Le ${\rm qe}$-module $\res L$ est la somme
orthogonale de ses composantes $p$-primaires, et pour tout premier
$p$, le morphisme \'evident $L^\sharp \rightarrow L_p^\sharp$ induit
une identification de la composante $p$-primaire de $\res L$
\`a $\res L_p$.

Notons que $\res L$ est anisotrope si et seulement si $\res L_p$ est
anisotrope pour tout premier $p$.  En effet, un \'el\'ement de
$\QQ/\ZZ$ est nul si et seulement si toutes ses composantes
$p$-primaires sont nulles.

Supposons de plus $E$ muni d'une structure de $\HH$-espace vectoriel
\`a droite telle que tout \'el\'ement $a$ de $\HH$ agisse sur $E$ par
une similitude orthogonale de rapport $\n(a)$.  On a alors $(x a)
\cdot y = x \cdot (y \overline{a})$ pour tous les $x,y$ dans $E$ et
$a$ dans $\HH$.  Si $L$ est stable par $\OOO$, il en va de m\^eme de
$L^\sharp$ car $\OOO$ est stable par la conjugaison $x \mapsto \overline{x}$, et $\res L$
est donc muni d'une structure de $\OOO$-module \`a droite. De m\^eme,
$L_p^\sharp$ et $\res L_p$ sont des $\OOO$-modules \`a droite, et
l'application naturelle $\res L \rightarrow \res L_p$ est
$\OOO$-lin\'eaire, pour tout premier $p$.

\medskip
Nous renvoyons \`a \cite{Vigneras80} pour les faits tr\`es classiques
suivants sur les alg\`ebres de quaternions. Soit $p$ un nombre premier
divisant $D_A$.  Alors $\OOO_p$ est l'unique ordre maximal de la
$\QQ_p$-alg\`ebre de quaternions $A \otimes_\QQ \QQ_p$ (une alg\`ebre
\`a division). De plus, la norme r\'eduite $\n : \OOO_p \rightarrow
\ZZ_p$ est surjective, et l'on peut donc choisir un \'el\'ement
$\pi_p$ de $\OOO_p$ v\'erifiant $\n (\pi_p)=p$.  Tout id\'eal (\`a
droite ou \`a gauche) de $\OOO_p$ est bilat\`ere, de la forme $\pi_p^n
\OOO_p$ pour un unique entier $n\geq 0$ (d'indice $p^{2n}$).  En
particulier, on a $p\OOO_p = \pi_p^2 \OOO_p$, $\pi_p \OOO_p$ est
l'id\'eal maximal de $\OOO_p$, et $\FF_{p^2} = \OOO_p /\pi_p\OOO_p$
est un corps d'ordre $p^2$. Nous noterons respectivement ${\rm
  Tr}_{\FF_{p^2}/\FF_p}$ et ${\rm N}_{\FF_{p^2}/\FF_{p}}$ la trace et
la norme de l'extension $\FF_{p^2}$ de $\FF_{p}=\ZZ_p/p\ZZ_p$.

D'apr\`es ces rappels, il existe un unique id\'eal \`a droite $\M$ de
$\OOO$ d'indice $D_A^2$~: c'est l'id\'eal v\'erifiant $\M_p=\pi_p
\OOO_p$ pour $p$ divisant $D_A$, et $\M_p = \OOO_p$ sinon.  C'est un
id\'eal bilat\`ere de $\OOO$ car $\M_p$ est bilat\`ere dans $\OOO_p$
pour tout $p$. Puisque sa norme est \'egale \`a $D_A$, l'id\'eal $\M$ est la
diff\'erente de $\OOO$, donc \'egal \`a $(\OOO^\sharp)^{-1}$ o\`u $I^{-1}=
\{x\in A\;:\; IxI\subset I\}$ pour tout $\ZZ$-r\'eseau $I$ de $A$.

\medskip
Puisque $\tr (\overline{x}\,x)=2\n (x)$ pour tout $x$ dans $\HH$,
l'ordre maximal $\OOO$ est un $\ZZ$-r\'eseau entier et pair de l'espace
euclidien $\HH$, de d\'eterminant $D_A^{\;\;2}$. Posons
$\N=\frac{1}{\sqrt{D_A}}\M\subset \HH$. C'est un sous-$\OOO$-module
bilat\`ere de $\HH$ v\'erifiant trivialement $\n(xa)=\n(ax)=\n(a)\n(x)$ pour tout $x \in \N$ et tout $a \in \OOO$.

\blemm\label{lem:residuordremax} 
\begin{enumerate} 
\item[(1)] Pour tout premier $p$ divisant $D_A$, le r\'esidu de
  $\OOO_p$ est isomorphe au groupe additif $\FF_{p^2}$ muni de la
  forme quadratique $x \mapsto \frac{1}{p} {\rm N}
  _{\FF_{p^2}/\FF_p}(x) \bmod \ZZ$.
\item[(2)] Le $\ZZ$-r\'eseau euclidien $\N$ est entier et pair, de dual
  $\N^\sharp = \frac{1}{\sqrt{D_A}} \OOO$, et pour tout premier $p$,
  il existe une isom\'etrie $\OOO_p$-lin\'eaire (\`a droite) entre
  $\N_p$ et $\OOO_p$.
\end{enumerate}
\elemm

\dem Soit $p$ un nombre premier divisant $D_A$.  La conjugaison $a
\mapsto \overline{a}$ de $\OOO$ induit une anti-involution de
$\OOO_p$, pr\'eservant n\'ecessairement son id\'eal maximal $\M_p$,
ainsi donc qu'un automorphisme $\FF_p$-lin\'eaire de $\FF_{p^2}$,
n\'ecessairement non trivial \`a cause de l'identit\'e $\n(a) = a
\overline{a}$, et donc \'egal \`a l'automorphisme de Frobenius $y\mapsto
y^p$.  Pour tout $x$ dans $\OOO_p$ d'image $y$ dans $\FF_{p^2}$, et
puisque $\pi_p\OOO_p\cap\ZZ_p=p\ZZ_p$, on a donc

\begin{equation}\label{normtrmodp}
\tr(x)\equiv {\rm Tr}_{\FF_{p^2}/\FF_p}(y) \bmod p\ZZ_p 
\hspace{5mm} \text{et}\hspace{5mm} 
\n(x) \equiv{\rm N}_{\FF_{p^2}/\FF_{p}}(y)  \bmod p\ZZ_p\;.
\end{equation}
De la congruence \eqref{normtrmodp} portant sur $\tr$,
on d\'eduit ais\'ement les relations bien connues
\begin{equation}\label{dieselloc}
  \tr \M_p \subset p\ZZ_p, \hspace{5 mm}
  \OOO_p^\sharp = p^{-1} \M_p = \pi_p^{-1}\OOO_p
  \hspace{5 mm}\text{et}\hspace{5 mm}
  \M_p^\sharp = p^{-1}\, \OOO_p\;.
\end{equation} 
La congruence \eqref{normtrmodp} portant sur $\n$, la relation
$\OOO_p^\sharp = \pi_p^{-1} \OOO_p$ ci-dessus, l'\'egalit\'e
$\n(\pi_p)=p$, et la multiplicativit\'e de la norme, entra\^inent
l'assertion (1) du lemme \ref{lem:residuordremax} .

Montrons l'assertion (2). Pour tout $x$ dans $\M$, l'\'el\'ement
$\frac{1}{D_A} \n(x)$ est dans $\ZZ_p$ pour tout premier $p$, et
donc dans $\ZZ$ : les r\'eseaux $\M$ et $\N$ sont entiers et pairs.
On a $\M^\sharp = D_A^{-1} \OOO$ d'apr\`es \eqref{dieselloc}, et donc
$\N^\sharp =\sqrt{D_A} \,\M^\sharp= \frac{1}{\sqrt{D_A}} \OOO$.  Soit
$p$ un nombre premier, il ne reste qu'\`a montrer qu'il existe une
bijection $\OOO_p$-lin\'eaire \`a droite $f_p : \OOO_p \rightarrow
\M_p$ v\'erifiant $\frac{1}{D_A}\n(f_p(x))\,=\n(x)$ pour
tout $x$ dans $\OOO_p$.  Il suffit de prendre pour $f_p$ la
multiplication \`a gauche par un \'el\'ement $x_p \in \OOO_p$ avec
$\n(x_p)=D_A$.  Un tel \'el\'ement existe par la surjectivit\'e
de $\n : \OOO_p \rightarrow \ZZ_p$.
\cqfd

\bprop \label{prop:existEhuit}
Il existe un $\OOO$-r\'eseau de $\HH \times \HH$ contenant
$\OOO \times \N$ et qui est isom\'etrique \`a ${\rm E}_8$ en tant que
$\ZZ$-r\'eseau euclidien. 
\eprop

\dem Notons $L$ le $\ZZ$-r\'eseau $\OOO\times \N$ de l'espace
euclidien $\HH\times\HH$. Il est stable par multiplication \`a droite
par $\OOO$, mais n'est pas n\'ecessairement un $\OOO$-r\'eseau car
$\N$ n'est pas libre de rang $1$ sur $\OOO$ en
g\'en\'eral.\footnote{Voir la partie \ref{subsect:ordmaxdiffprinc}
  pour des listes et des caract\'erisations de quand $\M$ (et donc $\N$)
  est libre de rang $1$ sur $\OOO$.}  Notons que $\ZZ$-r\'eseau
$L^\sharp$, et donc le groupe ab\'elien $L^\sharp/L$ de $\res L$, ont
des structures de $\OOO$-modules (\`a droite) telles que la projection
canonique $\operatorname{pr}: L^\sharp\ra L^\sharp/L$ soit un
morphisme de $\OOO$-module. Montrons qu'il existe un
sous-$\OOO$-module isotrope $I$ de $\res L$, d'ordre ${D_A}^2$ et tel
que le sous-$\OOO$-module $\Lambda=\operatorname{pr}^{-1}(I)$ de
$\HH\times\HH$ soit un $\OOO$-r\'eseau de $\HH\times\HH$.  Alors
$\Lambda$ est un $\ZZ$-r\'eseau euclidien entier et pair en dimension
$8$, de covolume $\covol \Lambda = \frac{\covol L} {[\Lambda:L]}=
\frac{{D_A}^2}{{D_A}^2}=1$ donc unimodulaire. Par unicit\'e, il est
isom\'etrique \`a ${\rm E}_8$, ce qui conclut.

Nous avons $\res L\, = \, \res \OOO \times \res
\N$ et $\res L_p\, = \, \res \OOO_p \times \res
\N_p$ pour tout premier $p$.  Il suffit donc de d\'efinir la
composante $p$-primaire $I_p$ de $I$ pour $p$ divisant $D_A$.  Fixons
un tel $p$ et identifions $\res \OOO_p$ et $\res \N_p$
au $\OOO$-module $\FF_{p^2}$ muni de la forme quadratique $x \mapsto
\frac{1}{p} {\rm N}_{\FF_{p^2}/\FF_p}(x) \bmod \ZZ$, ce qui est
loisible d'apr\`es les points (a) et (b) du lemme
\ref{lem:residuordremax}.  Soit $a_p\in\FF_{p^2}$ tel que ${\rm
  N}_{\FF_{p^2}/\FF_{p}}(a_p)=-1$, qui existe par la surjectivit\'e de
la norme pour les corps finis. Posons
\begin{equation}\label{eq:Ip}
I_p =\{(x_p,y_p)\in \FF_{p^2}\times\FF_{p^2}\;:\; y_p=
a_px_p\}\;.
\end{equation}
C'est un sous-$\OOO$-module de $\res L_p$ d'ordre $p^2$. Il est
isotrope car pour $(x_p, y_p) \in I_p$ on a
$$
{\rm N}_{\FF_{p^2}/\FF_{p}}(x_p)+ {\rm N}_{\FF_{p^2}/\FF_{p}}(y_p)
= {\rm N}_{\FF_{p^2}/\FF_{p}}(x_p) (1+ {\rm
  N}_{\FF_{p^2}/\FF_{p}}(a_p)) = 0\;.
$$
Enfin, $\Lambda =\operatorname{pr}^{-1}(I)$ est dans une suite exacte de
$\OOO$-modules
$$
0\longrightarrow \OOO \longrightarrow\Lambda
\longrightarrow \N^\sharp \longrightarrow 0\;,
$$
o\`u l'application $\Lambda \ra \N^\sharp$ est la restriction \`a
$\Lambda$ de la seconde projection $\OOO^\sharp\times \N^\sharp \ra
\N^\sharp$. En effet, cette application est surjective de noyau $\OOO
\times \{0\}$ car la seconde projection de $I_p$ dans $\res \N_p$
est bijective.  Mais on a $\N^\sharp = \frac{1}{\sqrt{D_A}} \OOO$
d'apr\`es le lemme \ref{lem:residuordremax} (b).  Donc $\N^\sharp$ est
libre de rang $1$ sur $\OOO$, et $\Lambda$ est un $\OOO$-r\'eseau, ce
qui conclut la d\'emonstration de la proposition.
\cqfd
\medskip

La proposition \ref{prop:E8}, et
donc le th\'eor\`eme \ref{theo:appBminbinhamform}, en d\'ecoulent.

\medskip

Notons $\udh$ le groupe unitaire du $\HH$-espace vectoriel
\`a droite $\HH \times \HH$ muni de la forme hamiltonienne
$f_0(u,v)=\n(u) + \n(v)$.  Ce groupe agit naturellement sur 
l'ensemble des $\OOO$-r\'eseaux de $\HH \times \HH$. 
\`A tout tel r\'eseau $\Lambda$, disons de base $(e_1,e_2)$, 
on associe la forme hamiltonienne binaire $(u,v) \mapsto f_0(e_1 u+e_2 v)$.
La classe de $\OOO$-\'equivalence de cette forme ne d\'epend 
que de $\Lambda$, et nous la notons $f_\Lambda$. 
Il d\'ecoule imm\'ediatement des d\'efinitions que pour deux r\'eseaux $\Lambda$ et $\Lambda'$, on a
l'\'egalit\'e $f_\Lambda = f_{\Lambda'}$ si, et seulement si, 
il existe $g \in \udh$ v\'erifiant $g(\Lambda)=\Lambda'$.
Comme toute forme hamiltonienne binaire d\'efinie positive est de la forme 
$f_0(e_1u+e_2 v)$ pour une base $(e_1,e_2)$ bien choisie du $\HH$-espace vectoriel $\HH \times \HH$, 
on en d\'eduit le r\'esultat suivant, analogue naturel d'un \'enonc\'e classique sur les r\'eseaux euclidiens.
(L'assertion portant sur le covolume a d\'ej\`a \'et\'e vue plus haut.)

\bprop L'application $\Lambda \mapsto f_\Lambda$ induit une bijection entre 
l'ensembles des $\udh$-orbites de $\OOO$-r\'eseaux dans $\HH \times \HH$
et l'ensemble des classes de $\OOO$-\'equivalence 
de formes hamiltoniennes binaires d\'efinies positives. Dans cette bijection, le discriminant de $f_\Lambda$ et le covolume de $\Lambda$
sont li\'es par la formule  $\sqrt{{\rm Covol}\, \Lambda} \,=\, -D_A \,\Delta(f_\Lambda)$. 
\cqfd
\eprop

\medskip

Notons $\mathcal{X}(\OOO)$ l'ensemble des $\OOO$-r\'eseaux de $\HH
\times \HH$ isom\'etriques \`a ${\rm E}_8$ en tant que r\'eseau
euclidien. Il est stable sous l'action de $\udh$.
\medskip

\bcoro 
\label{coro:udhbin}
\begin{enumerate}
\item L'ensemble $\mathcal{X}(\OOO)$ est non vide. 
\item L'application $f \mapsto f_\Lambda$ induit une bijection
entre $\udh \backslash \mathcal{X}(\OOO)$ et l'ensemble des classes de $\OOO$-\'equivalence 
de formes hamiltoniennes binaires d\'efinies positives de discriminant $-1/D_A$ qui r\'ealisent la (borne sup\'erieure d\'efinissant la) constante d'Hermite
$\gamma_2(\OOO)$.
\end{enumerate}
\ecoro

\dem
La premi\`ere assertion est la proposition  \ref{prop:E8}.
La seconde r\'esulte de l'analyse pr\'ec\'edant cette proposition et 
du th\'eor\`eme d\'ej\`a cit\'e de Vetchinkin \cite[Theo.~2]{Vetchinkin80},
\cqfd

\medskip

Bien que nous ne l'utiliserons pas, mentionnons que pour des raisons g\'en\'erales il n'y a qu'un
nombre fini de $\udh$-orbites dans $\mathcal{X}(\OOO)$.

\section{\'Etude du cas d'\'egalit\'e}
\label{sec:casegal}

Dans cette partie, nous explicitons d'abord certains \'el\'ements de 
$\mathcal{X}(\OOO)$ quand la diff\'erente de $\OOO$
est principale, ainsi donc que des formes hamiltoniennes binaires d\'efinies positives
r\'ealisant $\ga_2(\OOO)$. Ensuite, nous d\'eterminerons $\udh \backslash \mathcal{X}(\OOO)$ quand $\OOO$ est principal.

\medskip

Supposons donc que $\OOO$ est un ordre maximal de diff\'erente $\M$ principale.
Choisissons $\pi\in\OOO$ v\'erifiant
$\M=\pi\OOO=\OOO\pi$.  En particulier, nous avons $\n(\pi)=D_A$ et $\OOO^\sharp = \pi^{-1} \OOO$
(observer par exemple $\pi_p\in\pi\OOO_p^\times$ pour tout premier $p$ divisant $D_A$ et utiliser la formule \eqref{dieselloc}).
Notons $\E_\OOO$ l'ensemble des
$\lambda\in\OOO$ tels que $\n(\lambda)\equiv -1 \mod D_A$.
Cet ensemble est non vide par la surjectivit\'e de la norme
pour les corps finis et la formule \eqref{normtrmodp}. 
Pour $\lambda \in \E_\OOO$, posons
\begin{equation}\label{eq:defiLambdalambda}
\Lambda_{\lambda}=\{(\pi^{-1}u,\pi^{-1}v)\in\pi^{-1}\OOO\times \pi^{-1}\OOO
\;:\; \lambda u \equiv v\mod \pi\OOO\}\;.
\end{equation}
C'est un $\ZZ$-r\'eseau de $\HH \times \HH$ stable par $\OOO$ et contenant $\OOO \times \OOO$. 
L'anneau $\OOO/\pi\OOO$ \'etant commutatif (il est isomorphe \`a $\prod_{p\, |\, D_A}\FF_{p^2}$), 
 le r\'eseau $\Lambda_\lambda$ ne change pas si l'on remplace $\pi$ par $a \pi$, avec $a \in \OOO^\times$,
 dans sa d\'efinition. Autrement dit, $\Lambda_\lambda$ ne d\'epend pas 
 du choix du g\'en\'erateur $\pi$ de $\M$, ce qui justifie sa notation.
Enfin, il ne d\'epend que de la classe de $\lambda$ dans $\OOO/\pi\OOO$,
de sorte qu'il y a \'egalement un sens \`a d\'efinir $\Lambda_{\lambda}$
pour tout $\lambda \in \OOO/\pi\OOO$ avec $\n(\lambda) \equiv -1 \bmod D_A$.
En consid\'erant l'\'el\'ement $(\pi^{-1},\pi^{-1}\lambda)$ de $\Lambda_\lambda$, on constate l'\'equivalence $\Lambda_{\lambda} = \Lambda_{\lambda'}$ 
$\Leftrightarrow$ $\lambda \equiv \lambda' \bmod \pi\OOO$.
\par
\medskip

\bprop \label{prop:realgadeO} 
\label{prop:Lambdalambda} 
Supposons la diff\'erente de $\OOO$ principale et engendr\'ee par $\pi$. Soit $\lambda \in \E_\OOO$.\begin{enumerate} 
\item[(1)] Le $\ZZ$-r\'eseau $\Lambda_{ \lambda}$ appartient \`a $\mathcal{X}(\OOO)$. 
\item[(2)] L'application de $\HH\times\HH$ dans
$\RR$ d\'efinie par
$$ 
f_{\pi; \lambda} :(u,v)\mapsto  \frac{\n(\lambda)+1}{D_A}\,\n(u)+ \tr(\,\overline{u}\,\pi^{-1}\lambda\,v)
+\,\n(v)
$$ 
est une forme hamiltonienne binaire d\'efinie positive, qui r\'ealise la
borne sup\'erieure d\'efinissant $\ga_2(\OOO)$. 
\end{enumerate}
Enfin, tout \'el\'ement de $\mathcal{X}(\OOO)$ contenant $\OOO \times \OOO$ est de la forme $\Lambda_\lambda$ pour $\lambda \in \E_\OOO$.
\eprop
\medskip

\medskip

\dem Le $\ZZ$-r\'eseau $\Lambda_{ \lambda}$ contient $\OOO \times \OOO$ avec indice $|\OOO/\pi\OOO|=D_A^2 = {\rm Covol} (\OOO \times \OOO)$, il est donc 
bien de covolume $1$. Pour $(\pi^{-1}u,\pi^{-1}v) \in \Lambda_\lambda$, le rationnel 
$$\n(\pi^{-1}u)+\n(\pi^{-1}v)=\n(\pi)^{-1} (\n(u)+\n(v))=\frac{1}{D_A}(\n(u)+\n(v))$$ est un entier, 
\`a cause de la congruence  $\n(u)+\n(v) \equiv \n(u) +\n(\lambda) \n(u) \equiv 0 \bmod D_A$ pour $\lambda \in \E_\OOO$.
Cela montre que $\Lambda_{ \lambda}$ est entier pair. C'est trivialement un sous-$\OOO$-module (\`a droite) de $\HH \times \HH$, 
dont les \'el\'ements  $e_1=(\pi^{-1},\pi^{-1}\lambda)$ et $e_2=(0,1)$ 
constituent une $\OOO$-base : c'est un $\OOO$-r\'eseau, et nous avons montr\'e le (1).
On constate les \'egalit\'es $f_0(e_1) =  \frac{1+\n(\lambda)}{D_A}$, $f_0(e_2) =1$ et 
$h(e_1,e_2)= \pi^{-1} \lambda$. Pour $u,v \in \HH$, nous avons donc 
$f_0(e_1 u + e_2 v) = f_{\pi;\lambda}(u,v)$. L'assertion (2) d\'ecoule alors de l'assertion (1)
et du point (2) du corollaire \ref{coro:udhbin}.

Montrons la derni\`ere assertion. 
On raisonne comme dans la d\'emonstration de la proposition \ref{prop:existEhuit}, 
en rempla\c{c}ant $\OOO \times \mathcal{N}$ par $\OOO \times \OOO$.
Nous avons ${\rm res}\, \OOO = \pi^{-1} \OOO/\OOO$ et l'anneau $\OOO/\pi\OOO$ est 
produit direct sur les premiers $p$ divisant $D_A$ des corps finis $\OOO_p/\pi\OOO_p=\FF_{p^2}$.
Un sous-$\OOO/\pi\OOO$-module $I$ totalement isotrope de ${\rm res}\, \OOO \times {\rm res}\, \OOO$ 
intersecte trivialement le sous-espace anisotrope $\{0\} \times {\rm res}\, \OOO_p$ pour $p$ divisant $D_A$.
Si $I$ est en outre libre de rang $1$ sur $\OOO/\pi\OOO$ (ou ce qui revient au m\^eme, de cardinal $|\OOO/\pi\OOO|$), 
il est engendr\'e sur $\OOO/\pi\OOO$ par la classe d'un \'el\'ement $e$ de $\pi^{-1}\OOO \times \pi^{-1}\OOO$, disons   
$e=(\pi^{-1}\mu',\pi^{-1}\mu)$, avec $\mu',\mu \in \OOO$ v\'erifiant $\n(\mu')+\n(\mu) \equiv 0 \bmod D_A$
et $\mu' \in \OOO_p^\times$ pour tout $p$ divisant $D_A$. On en d\'eduit $\mu'\in \OOO^\times$ :  on peut donc
supposer $\mu'=1$. Par d\'efinition, l'image inverse de $I$ par la projection canonique 
$\pi^{-1} \OOO \times \pi^{-1} \OOO \rightarrow {\rm res} \,\OOO \times {\rm res}\,\OOO$ est $e \OOO + \OOO \times \OOO=\Lambda_\mu$.
\cqfd

\noindent
\brema\label{rem:constructE8}{\rm Supposons $D_A=2$ et $\OOO=\ZZ\frac{1+i+j+k}{2}+i\ZZ+j\ZZ+k\ZZ$ ({\it ordre
de Hurwitz}).  Nous pouvons prendre $\pi=i+1$ (car $\n(1+i)=2$) et
$\lambda =1$ (car $1=-1\mod 2$). Nous retrouvons alors la r\'ealisation quaternionique usuelle $\{(u',v')\in
\OOO\times \OOO \;:\; u'\equiv v'\mod (1+i)\OOO\}$ du r\'eseau euclidien
${\rm E}_8$ : voir
\cite[Prop.~8.2.2]{Martinet03} (cette construction diff\`ere de la n\^otre d'une homoth\'etie car 
elle utilise le produit scalaire $\frac{1}{2}
\tr(\overline{x}y)$ sur $\HH$). }\par
\erema

\bigskip
D\'eterminons maintenant
$\udh \backslash \mathcal{X}(\OOO)$ lorsque $\OOO$ est principal, ce
qui se produit si et seulement si $h_A=1$, ou de mani\`ere \'equivalente
si $D_A = 2, 3, 5, 7$ ou $13$. En particulier, $D_A$ est un nombre
premier, que nous noterons simplement $p$. La
proposition \ref{prop:nonuniintro} de l'introduction d\'ecoule du
r\'esultat suivant.

\bprop \label{propcaseg} Pour $D_A \leq 7$, il
existe une unique $\udh$-orbite de $\OOO$-r\'eseaux de $\HH \times
\HH$ isom\'etriques \`a ${\rm E}_8$. Pour $D_A=13$, il existe exactement
deux telles orbites.  
\eprop

Comme $\OOO$ est principal, sa diff\'erente l'est aussi, et on en fixe comme pr\'ec\'edemment
un g\'en\'erateur $\pi$. Rappelons que si $h_A=1$, la formule de masse d'Eichler
\cite[p.~103]{Eichler38} donne
$$
\frac{1}{|\OOO^\times|}=\frac{p-1}{24}\;.
$$ 

\'Ecrivons l'entier $24$ sous la forme $(p-1) p^n m$ avec $m$ un
entier premier \`a $p$.  Le morphisme d'anneaux $\OOO_p \rightarrow
\OOO_p/\pi \OOO_p = \FF_{p^2}$ induit un morphisme de groupes
$\OOO_p^\times \rightarrow \FF_{p^2}^\times$ de noyau $1+\pi
\OOO_p$.  Ce dernier est un pro-$p$-groupe et $\FF_{p^2}^\times$ est
d'ordre premier \`a $p$.  En composant $\OOO_p^\times \rightarrow
\FF_{p^2}^\times$ et l'inclusion canonique $\OOO^\times \rightarrow
\OOO_p^\times$, nous obtenons un morphisme de groupes
$$
\varepsilon: \OOO^\times \rightarrow \FF_{p^2}^\times\;,
$$
dont l'image est incluse dans le sous-groupe d'ordre $p+1$ des
\'el\'ements de norme $1$ de $ \FF_{p^2}^\times$.  Le lemme suivant en
d\'ecoule.

\blemm\label{lem:surjeps} 
Nous avons $|{\rm ker}\, \varepsilon|\,=\,p^n$, $|{\rm Im} \,\varepsilon|
\,=\, m$, et $m$ divise $p+1$. \cqfd 
\elemm

D'apr\`es le dernier point de la proposition
\ref{prop:realgadeO}, et puisqu'ici $D_A=p$ est premier, les
r\'eseaux $\Lambda\in\mathcal{X}(\OOO)$ contenant $\OOO \times \OOO$ sont 
ceux de la forme $\Lambda_{x}$ pour un unique $x$ dans 
$\OOO/\pi\OOO = \mathbb{F}_{p^2}$ v\'erifiant 
$$
{\rm N}_{\FF_{p^2}/\FF_p}(x) = x^{1+p} = -1.
$$
Notons $X_p \subset \FF_{p^2}^\times$ le sous-ensemble des $x$ ci-dessus.
Nous avons $|X_p|=p+1$. 
On munit $X_p$ d'une structure de $\OOO^\times$-ensemble par la formule $(a,x) \mapsto \varepsilon(\pi a \pi^{-1}) x$. 
Cette formule a un sens car la relation $\pi\OOO = \OOO \pi$ montre que $x \mapsto \pi x \pi^{-1}$ est un 
automorphisme de l'anneau $\OOO$, et en particulier du groupe $\OOO^\times$. 
Observons que le $\OOO^\times$-ensemble $X_p$ ainsi d\'efini ne d\'epend pas du choix du g\'en\'erateur $\pi$ de la diff\'erente de $\OOO$, car 
l'anneau $\OOO/\pi\OOO =\FF_{p^2}$ est commutatif.

Pour $x \in X_p$, le $\ZZ$-r\'eseau 
$\Lambda_{x}$ contient l'\'el\'ement $\alpha_0 = (1,0)$ qui v\'erifie
$\alpha_0 \cdot \alpha_0=2$. Pour tout $a$ dans $\OOO^\times$, notons ${\rm m}_a$ l'\'el\'ement de
$\udh$ d\'efini par ${\rm m}_a(x,y)=(x,ay)$.  Le morphisme $a \mapsto
{\rm m}_a$ identifie $\OOO^\times$ au sous-groupe de $\udh$ fixant
$\alpha_0$ et pr\'eservant $\OOO \times \OOO$. Pour $a \in \OOO^\times$, nous
avons l'identit\'e
${\rm m}_a(\pi^{-1},\pi^{-1}\lambda)=(\pi^{-1},\pi^{-1} \pi a \pi^{-1} \lambda)$. 
La d\'efinition \eqref{eq:defiLambdalambda} montre donc
$${\rm m}_a(\Lambda_{x}) =\Lambda_{\varepsilon(\pi a \pi^{-1}) x}$$
pour $x \in X_p$ et $a \in \OOO^\times$.
Notons enfin $\Y (\OOO)$ l'ensemble des couples
$(\alpha,\Lambda)$ avec $\Lambda$ dans $\X(\OOO)$, $\alpha
\in \Lambda$, et $\alpha \cdot \alpha =2$, muni de l'action diagonale
de $\udh$.  Les observations ci-dessus montrent que l'application $X_p \rightarrow
\Y(\OOO)$, d\'efinie par $x \mapsto (\alpha_0,\Lambda_{x}),$ et le
morphisme $\OOO^\times \rightarrow \udh$, d\'efini par $a \mapsto {\rm
  m}_a$, d\'efinissent un morphisme de groupo\"{\i}des entre le
$\OOO^\times$-ensemble $X_p$ et le $\udh$-ensemble $\Y(\OOO)$.

\blemm\label{lem:orbalphaL} 
Le morphisme de groupo\"{\i}des $X_p \rightarrow \Y(\OOO)$ ci-dessus est
une \'equivalence (de cat\'egorie).  
\elemm

\dem Il y a deux points \`a d\'emontrer : (i) pour tout $y$ dans
$\Y(\OOO)$, il existe $u$ dans $\udh$ et $x$ dans $X_p$ tels que
$uy = (\alpha_0,\Lambda_{x})$; (ii) pour tout $x$ dans $X_p$, le
stabilisateur de $(\alpha_0,\Lambda_{x})$ dans $\udh$ est le sous-groupe
des ${\rm m}_a$ avec $a \in \OOO^\times$ et $\varepsilon(\pi a \pi^{-1})=1$.

Soient $\Lambda\in\X(\OOO)$, $\alpha \in \Lambda-\{0\}$ et $N =
\alpha \OOO$. Pour $x$ dans $\OOO$, la relation $\alpha x\, \cdot\, \alpha x \,=\, \n(x)\,
\alpha \cdot \alpha \,=\,\mu\,(x \cdot x)$, avec $\mu = \frac{\alpha\cdot \alpha}{2}=f_0(\alpha)$, montre que le
$\ZZ$-r\'eseau euclidien de rang $4$ sous-jacent \`a $N$ est
isom\'etrique \`a $\sqrt{\mu}\OOO$, et donc de covolume $\mu^2 \covol
\OOO$. Supposons maintenant $\alpha \cdot \alpha = 2$.  Dans ce cas,
$N$ est isom\'etrique \`a $\OOO$.  Le r\'esidu de $\OOO$ \'etant
anisotrope, le seul r\'eseau entier pair de $N \otimes_\ZZ \QQ$ contenant
$N$ est $N$ lui-m\^eme.  En particulier, le r\'eseau $N$ est satur\'e
dans $\Lambda$ : on a $\Lambda \cap (N \otimes_\ZZ \QQ) = N$ et le groupe
ab\'elien $\Lambda/N$ est sans $\ZZ$-torsion.  Comme ${\rm E}_8$ est
unimodulaire, on en d\'eduit que l'orthogonal $N^\perp$ de $N$ dans
$\Lambda$ est de m\^eme covolume que $N$
\cite[Prop.~B~2.2.~(d)]{CheLan19}.  Mais $N^\perp$ est stable par
$\OOO$ car $N$ l'est.  Il est donc libre de rang $1$ sur $\OOO$ car
$\OOO$ est principal.  Nous pouvons donc \'ecrire $N^\perp = \beta
\OOO$ avec $\beta \in \Lambda$. L'\'egalit\'e des covolumes de $N$
et $N^\perp$ implique $\beta \cdot \beta =2$.  Ainsi, le couple
$(\alpha,\beta)$ est une $\HH$-base orthonorm\'ee de $\HH\times\HH$,
et donc quitte \`a remplacer $\Lambda$ par $u (\Lambda)$ o\`u $u$ est
l'unique \'el\'ement de $\udh$ envoyant $\alpha$ sur $\alpha_0=(1,0)$
et $\beta$ sur $(0,1)$, nous pouvons supposer que $\Lambda$, qui est
un $\OOO$-r\'eseau, contient $\OOO\times\OOO$, et que l'orthogonal de
$(1,0) \OOO$ dans $\Lambda$ est $\{0\} \times \OOO$.  En particulier, 
d'apr\`es la derni\`ere assertion de la proposition \ref{prop:realgadeO},
le r\'eseau $\Lambda$ est de la forme $\Lambda_x$ pour $x$ dans $X_p$~: le premier
point (i) en d\'ecoule.

De plus, le sous-groupe de $\udh$ fixant $\alpha_0$ et $\Lambda$
pr\'eserve l'orthogonal de $\alpha_0 \OOO$ dans $\Lambda$,
c'est-\`a-dire $\{0\} \times \OOO$~: c'est donc l'ensemble des ${\rm
  m}_a$ avec $a$ dans $\OOO^\times$.  La relation ${\rm m}_a (\Lambda_{x}) =
\Lambda_{\epsilon(\pi a \pi^{-1})x}$, et la propri\'et\'e $\Lambda_{x}=\Lambda_{y}$ $\Leftrightarrow$ $x=y$,
montrent que ${\rm m}_a (\Lambda_{x})= \Lambda_x$
\'equivaut \`a $\epsilon(\pi a \pi^{-1})=1$. Le second point (ii) en d\'ecoule.  
\cqfd

\bigskip
Pour $\Lambda \subset \HH\times\HH$, notons $\operatorname{U}
(\Lambda)$ le stabilisateur de $\Lambda$ dans $\udh$ (le {\it groupe
  unitaire de $\Lambda$}) et ${\rm R}(\Lambda)$ l'ensemble des $\alpha
\in \Lambda$ v\'erifiant $\alpha \cdot \alpha = 2$ (les {\it racines}
de $\Lambda$).  Nous avons $|{\rm R}(\Lambda)|=240$ pour $\Lambda$
dans $\mathcal{X}(\OOO)$, car ${\rm E}_8$ admet $240$ racines.  Le
groupe $\operatorname{U}(\Lambda)$ agit naturellement sur ${\rm
  R}(\Lambda)$, et on note ${\rm r}_\Lambda$ le nombre d'orbites pour
cette action.  Ce nombre ne d\'epend que de la $\udh$-orbite de
$\Lambda$.

\blemm \label{lem:tradorbalphaL} 
Le groupe $\udh$ admet $\frac{p+1}{m}$ orbites dans $\Y(\OOO)$, avec
stabilisateurs d'ordre $p^n$. Par cons\'equent, nous avons $\sum_{[\Lambda] \in \udh
  \backslash \X(\OOO)} {\rm r}_\Lambda = \frac{p+1}{m}$, et $240 = {\rm r}_\Lambda \frac{|\operatorname{U} 
(\Lambda)|}{p^n}$ pour tout $\Lambda$ dans $\mathcal{X}(\OOO)$.  
\elemm

\dem 
Le $\OOO^\times$-ensemble ${\rm X}_p$ a clairement $(p+1)/m$ orbites,
et des stabilisateurs isomorphes \`a ${\rm ker}\, \varepsilon$.
L'\'enonc\'e est donc une cons\'equence des lemmes \ref{lem:orbalphaL}
et \ref{lem:surjeps}.  
\cqfd

\medskip
\noindent
{\bf D\'emonstration de la proposition \ref{propcaseg}.}  
Pour $p=2,3,5$, nous avons $m=p+1$.  D'apr\`es le lemme
\ref{lem:tradorbalphaL}, le groupe $\udh$ agit transitivement sur
$\mathcal{X}(\OOO)$, et pour tout $\Lambda$ dans $\mathcal{X}(\OOO)$,
nous avons ${\rm r}_\Lambda=1$.

\brema Pour $p=2,3,5$ et $\Lambda$ dans $\mathcal{X}(\OOO)$ nous avons montr\'e
$|\operatorname{U}(\Lambda)| \,=\, 240 \,p^n$, soit
encore $$|\operatorname{U}(\Lambda)|\,\, =\,\, 1920,\,\, 720,\,\,
240 \hspace{.3 cm} \text{pour}\hspace{.3 cm}
p\,\,=\,\,2,\,3,\,5 \hspace{.3 cm}\text{respectivement.}$$ C'est un
point de d\'epart pour une \'etude plus fine du groupe
$\operatorname{U}(\Lambda)$ que nous laissons au lecteur.  Par
exemple pour $p=3$ (resp. $p=5$), on peut montrer que
$\operatorname{U}(\Lambda)$ est une extension centrale du groupe
altern\'e $\mathfrak{A}_6$ (resp. du groupe sym\'etrique
$\mathfrak{S}_5$) par $\ZZ/2\ZZ$. \footnote{On pourra d'abord observer que le groupe unitaire de
$\Lambda \otimes \ZZ/2\ZZ$ est isomorphe \`a ${\rm Sp}_4(\ZZ/2\ZZ)
\simeq \mathfrak{S}_6$. De plus, les fibres non vides de la projection canonique 
${\rm R}(\Lambda) \rightarrow \Lambda \otimes \ZZ/2\ZZ$ 
sont de la forme $\{\alpha,-\alpha\}$ (une propri\'et\'e classique de ${\rm E}_8$). Ainsi, pour $\alpha \in {\rm R}(L)$, 
le sous-groupe de ${\rm U}(\Lambda)$ fixant $\alpha \bmod 2 \Lambda$ 
est d'orde $2p^n$ d'apr\`es le lemme \ref{lem:tradorbalphaL}. Le noyau du morphisme naturel de ${\rm U}(\Lambda)$ 
vers le groupe unitaire de $\Lambda \otimes \ZZ/2\ZZ$ est un $2$-groupe, il est donc \'egal \`a $\{\pm 1\}$. 
On conclut car un sous-groupe d'indice $2$ (resp. $6$) de $\mathfrak{S}_6$ est isomorphe \`a $\mathfrak{A}_6$ (resp. $\mathfrak{S}_5$).
}  
\erema

Supposons maintenant $p=7$. Alors $m=4=\frac{p+1}{2}$ et $n=0$.
Consid\'erons l'\'el\'ement $\tau$ de $\udh$ d\'efini par $\tau(x,y) =
(y,x)$.  Il pr\'eserve $\OOO \times \OOO$ et permute donc \'egalement
les $\Lambda_{x}$ pour $x$ dans $X_p$. 
Nous avons en fait \'evidemment
$\tau(\Lambda_{x}) = \Lambda_{x^{-1}}$. 

\noindent\begin{minipage}{8.9cm} 
Il suffit donc pour conclure de
remarquer (voir le dessin ci-joint) que si $G$ est le sous-groupe des
permutations de l'ensemble $X_7\subset \FF_{49}^\times$ engendr\'e par $x \mapsto
x^{-1}$ et les multiplications par un \'el\'ement d'ordre $4$, alors
$G$ agit transitivement sur $X_7$. 
\end{minipage} 
\begin{minipage}{6cm}
\begin{center}
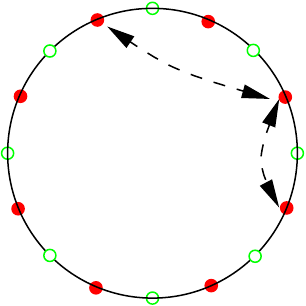
\end{center}
\end{minipage} 

\medskip
Consid\'erons enfin le cas plus d\'elicat $p=13$.  Alors $m=2$, $n=0$
et $\OOO^\times=\{\pm 1\}$.  Fixons $\Lambda$ dans $\X(\OOO)$.
D'apr\`es le lemme \ref{lem:tradorbalphaL}, nous avons
\begin{equation} \label{eq:classp13} 
240 = {\rm r}_\Lambda \;|\operatorname{U}(\Lambda)|\hspace{.5cm}
\text{et}\hspace{.5cm}
\sum_{ [\Lambda] \in \udh \backslash \mathcal{X}(\OOO)} {\rm r}_\Lambda =7\;.
\end{equation}
Si $\udh$ agit transitivement sur $\X(\OOO)$, nous avons ${\rm
  r}_\Lambda=7$ pour tout $\Lambda\in\X(\OOO)$ et donc $240 = 7
\;|\operatorname{U}(\Lambda)|$, une contradiction car $240$ n'est
pas divisible pas $7$. Donc $\udh$ admet au moins deux orbites
distinctes, puis ${\rm r}_\Lambda \leq 6$ et $40 \leq
|\operatorname{U}(\Lambda)|\leq 240$ pour tout $\Lambda\in\X(\OOO)$.

Soit $\lambda\in\OOO$ tel que $\n(\lambda) = 12$ ; l'arithm\'etique des
quaternions montre qu'il existe exactement $|\OOO^\times| \cdot 7
\cdot 4 = 56$ tels \'el\'ements $\lambda$.  Nous avons $\n(\lambda)
\equiv -1 \bmod p$, de sorte que $\lambda\in\E_\OOO$. 
La d\'emonstration de la proposition  \ref{prop:Lambdalambda} assure que dans la $\OOO$-base
de $\Lambda_\lambda$ d\'efinie par $e_1=(\pi^{-1},\pi^{-1}\lambda)$ et $e_2=(0,1)$, 
on a pour tous
$u,v \in \OOO$
$$f_0(e_1u + e_2 v) =  f_{\pi;\lambda}(u,v)\,=\,\n(u)\,+\,\tr(\,\overline{u}\,b\,v\,)\,+\,\n(v),
\hspace{.3cm}\text{avec}\hspace{.3cm} b\,=\, \pi^{-1} \lambda\;.
$$ 
En particulier, le cardinal de $\operatorname{U}(\Lambda_{\lambda})$ est
le nombre de couples de racines $(\alpha,\beta)$ de $\Lambda_{\lambda}$
de produit scalaire $\HH$-hermitien
$h(\alpha,\beta)=\pi^{-1}\lambda$. Il est facile d'\'enum\'erer ces couples \`a l'aide
d'un ordinateur~: il suffit d'\'enum\'erer ${\rm R}(\Lambda_{\lambda})$,
ou ce qui revient au m\^eme, les $(u,v)$ dans $\OOO \times \OOO$ tels
que $\n(u)+\n(v)=13$ et $\lambda u \equiv v \bmod \pi\OOO$.  D\'ecrivons le
r\'esultat. Nous pouvons prendre (voir par exemple
\cite[p.~98]{Vigneras80}) pour $A$ l'alg\`ebre de quaternions $\QQ \,+
\,\QQ \,i \,+\, \QQ\, j\, + \,\QQ \,k$ avec $i^2=-2$, $j^2=-13$ et
$ij=k=-ji$, et pour $\OOO$ l'ordre maximal $\ZZ \,+ \,\ZZ \,i\, +
\,\ZZ \,\frac{1+i-j}{2}\, +\,\ZZ \,\frac{2+i+k}{4}$.  
Nous pouvons alors prendre $\pi = j$.  
Consid\'erons les deux \'el\'ements de $\OOO$ de norme $12$
\begin{equation}
\label{lambda1lambda2}
\lambda_1=\frac{3+j-k}{2} \hspace{.5cm} \text{et}\hspace{.5cm}
\lambda_2=\frac{1+2i-j+k}{2}\;.
\end{equation}
L'ordinateur nous dit que $\Lambda_{\lambda_1}$ (resp.~$\Lambda_{\lambda_2}$) 
contient exactement $48$ (resp.~$120$) couples de racines $(\alpha,
\beta)$ de produit scalaire $\HH$-hermitien $j^{-1} \lambda_1$ (resp.
$j^{-1} \lambda_2$). Cela montre
$$
|\operatorname{U}(\Lambda_{\lambda_1})| = 48 \hspace{.5cm} 
\text{et}\hspace{.5cm}  |\operatorname{U}(\Lambda_{\lambda_2})| = 120\;.
$$ 
Il d\'ecoule de la formule \eqref{eq:classp13} les \'egalit\'es ${\rm
  r}_{\Lambda_{\lambda_1}} = 5$, ${\rm r}_{\Lambda_{\lambda_2}} = 2$,
puis ${\rm r}_{\Lambda_{\lambda_1}} + {\rm r}_{\Lambda_{\lambda_2}} =
7$, et donc $\udh\backslash \X(\OOO) = \{ \Lambda_{\lambda_1},
\Lambda_{\lambda_2}\}$. Ceci d\'emontre la proposition \ref{propcaseg}.
\cqfd

\medskip

D'apr\`es le corollaire \ref{coro:udhbin} et la proposition \ref{propcaseg}, nous avons d\'etermin\'e,
\`a $\OOO$-\'equivalence pr\`es, les formes hamiltoniennes binaires d\'efinies positives de d\'eterminant $-\frac{1}{D_A}$ (ou ce qui revient au m\^eme, de ``covolume $1$'') 
r\'ealisant $\ga_2(\OOO)$ lorsque $\OOO$ est principal. D\'ecrivons les formes trouv\'ees.
La principalit\'e de $\OOO$ enta\^ine que la norme r\'eduite $\n : \OOO \rightarrow \ZZ_{\geq 0}$ est surjective.
Le sous-ensemble $\E'_\OOO=\{ \lambda \in \OOO\;:\; \n(\lambda) = D_A-1\}$ de $\E_\OOO$ est donc non vide.
Pour $\lambda \in \E'_\OOO$, nous avons simplement $$f_{\pi;\lambda}(u,v) = \n(u) + \tr( \overline{u} \,\pi^{-1} \lambda \,v) + \n(v).$$ 
Pour $D_A=2,3,5$ et $7$, cette forme est donc l'unique forme de covolume $1$, \`a $\OOO$-\'equivalence pr\`es, r\'ealisant $\ga_2(\OOO)$.
Pour $D_A=13$, il y a deux classes d'\'equivalences, chacune \'etant repr\'esent\'ee par une telle forme. Explicitons des choix possibles de $\pi$ et $\lambda$ dans chacun des cas. \par\smallskip

\begin{itemize}
\item Pour $D_A=2$ et $\OOO$ l'ordre de Hurwitz usuel, les choix $\pi = 1+i$ et $\lambda=1$ conduisent \`a $\pi^{-1}\lambda = \frac{1-i}{2}$, conform\'ement \`a \cite{Speiser32}. \par\smallskip
\item Pour $D_A=3, 7$, la $\QQ$-alg\`ebre $A$ est engendr\'ee par deux \'el\'ements $i$ et $j$ v\'erifiant $i^2=-1$, $j^2=-D_A$ et $ij=-ji$. L'\'el\'ement $i$ normalise 
$\ZZ[\frac{1+j}{2}]$, de sorte que l'on peut choisir l'ordre maximal $\OOO$ contenant $i$ et $\frac{1+j}{2}$, et prendre $\pi=j$. On constate que $\lambda = 1+i$ (cas $D_A=3$) et 
$\lambda = 2+i \frac{1+j}{2}$ (cas $D_A=7$) conviennent, et conduisent respectivement \`a $\pi^{-1}\lambda = \frac{k-j}{3}$ et $\pi^{-1}\lambda = \frac{k-i-4j}{14}$, o\`u l'on a pos\'e $k=ij$.\par \smallskip
\item Pour $D_A=5,13$, la $\QQ$-alg\`ebre $A$ est engendr\'ee par deux \'el\'ements $i$ et $j$ v\'erifiant $i^2=-2$, $j^2=-D_A$ et $ij=-ji$. Dans le cas $D_A=5$, on peut prendre pour $\OOO$
tout ordre maximal contenant $j$, puis $\pi=j$ et $\lambda=2$, auquel cas nous avons $\pi^{-1}\lambda = - \frac{2 }{5} j$. Dans le cas $D_A=13$, l'analyse faite dans la d\'emonstration de la proposition
pr\'ec\'edente montre que, pour le choix de $\OOO$ de cette d\'emonstration et pour $\pi=j$, les deux classes sont donn\'ees par les \'el\'ements $\lambda=\lambda_i$, avec $i=1,2$, de la formule \eqref{lambda1lambda2}.

\end{itemize}

\brema \label{remxo}
Il serait int\'eressant de poursuive cette analyse en d\'eterminant 
des repr\'esentants de $\udh \backslash \mathcal{X}(\OOO)$ pour d'autres discriminants $D_A$, 
et aussi d'expliciter une formule donnant la masse 
du groupo\"ide $\udh \backslash \mathcal{X}(\OOO)$ pour un ordre maximal $\OOO$ g\'en\'eral. 
\erema

\section{Ordres maximaux et r\'eseaux euclidiens}
\label{subsect:ordreetreseaux}

Dans cette partie, nous rappelons certains aspects de la
correspondance classique (voir par exemple \cite{Latimer37,Peters69}
et \cite[Chap.~22]{Voight18}) entre classes de conjugaison d'ordres
maximaux de $A$ et certaines formes quadratiques ternaires. Nous
utiliserons cette correspondance dans la partie
\ref{subsect:ordmaxdiffprinc} pour donner des caract\'erisations des
ordres maximaux des alg\`ebres de quaternions rationnelles d\'efinies dont
la diff\'erente est principale.

Fixons un entier $d$ strictement positif sans facteur carr\'e, ayant
un nombre impair de facteurs premiers. Fixons aussi une alg\`ebre de
quaternions $A$ sur $\QQ$, qui est d\'efinie, de discriminant
$d$. Notons $\Rcal(d)$ l'ensemble des classes d'isom\'etrie de
r\'eseaux euclidiens entiers pairs $L$ de rang $3$, de d\'eterminant
$2d^2$, et tels que, pour tout premier $p$ impair divisant $d$, le
r\'esidu de $L_p$ est anisotrope.

Remarquons que si $L$ est dans $\mathcal{R}(d)$, et si $p$ est un
premier impair divisant $d$, nous avons $\res L_p\, \simeq \,\res
\OOO_p$, o\`u $\OOO$ d\'esigne un ordre maximal quelconque de $A$. En effet, le groupe ab\'elien $\res L_p$ est d'ordre $p^2$ et
ne peut \^etre isomorphe \`a $\ZZ/p^2\ZZ$, car sa $p$-torsion serait
constitu\'ee d'\'el\'ements isotropes. C'est donc un $\ZZ/p\ZZ$-espace
vectoriel de dimension $2$ muni d'une forme quadratique anisotrope,
n\'ecessairement \`a valeurs dans $(\frac{1}{p}\ZZ)/\ZZ \simeq \ZZ/p\ZZ$
car $p$ est impair.  Il n'y a qu'une telle forme \`a isom\'etrie
pr\`es, donn\'ee au lemme \ref{lem:residuordremax} (1).

\bprop 
\label{carunO}
L'application qui \`a un ordre maximal $\OOO$ de $A$ associe le
$\ZZ$-r\'eseau $L(\OOO)$ de ses \'el\'ements de trace nulle, muni du
produit scalaire $(x,y)\mapsto \tr(\overline{x}\,y)$, induit une
bijection de l'ensemble des classes de conjugaison d'ordres maximaux
de $A$ dans $\Rcal(d)$. En particulier,  $|\Rcal(d)|=t_A$.
\eprop

\dem Montrons tout d'abord que cette application est bien
d\'efinie. Soit $L=L(\OOO)$ Rappelons que $\tr:\OOO\ra\ZZ$ est
surjective puisque $\OOO$ est maximal : c'est \'evident si $\OOO_2
\simeq \M_2(\ZZ_2)$ et cela d\'ecoule de la surjectivit\'e de ${\rm
  Tr}_{\FF_4/\FF_2} : \FF_4 \rightarrow \FF_2$ et de la formule
\eqref{normtrmodp} concernant la trace si $2$ divise $d$.  Ainsi, $\ZZ
1+L$ est un $\ZZ$-r\'eseau d'indice $2$ dans $\OOO$. Son d\'eterminant
est donc $4\det \OOO=4d^2$ puisque $\OOO$ est maximal, et $L$ est bien
de d\'eterminant $2d^2$. Enfin, si $p$ est un premier impair divisant
$d$, nous avons $\res L_p\, \simeq \, \res(\ZZ \,1 \,\oplus L)_p\, =\,
\res \OOO_p,$ donc $\res L_p$ est anisotrope.

\medskip 
Montrons l'injectivit\'e. Soient $\OOO$ et $\OOO'$ deux ordres
maximaux de $A$ tels que les $\ZZ$-r\'eseaux euclidiens $L=L(\OOO)$ et
$L'= L(\OOO')$ soient isom\'etriques. Notons $A_0$ l'espace
quadratique des quaternions purs de $A$. Soit $u :A_0\ra A_0$ une
isom\'etrie telle que $u(L)=L'$. Par un r\'esultat classique (voir par
exemple \cite[p.~11 et 6]{Vigneras80}), $u$ est la conjugaison par un
\'el\'ement de $A-\{0\}$.  Quitte \`a remplacer $\OOO'$ par un
conjugu\'e, nous pouvons donc supposer $L=L'$. Les $\ZZ$-r\'eseaux
$\OOO$ et $\OOO'$ contiennent alors tous deux le $\ZZ$-r\'eseau $\ZZ
1+L$, avec l'indice $2$.  Cela montre $\OOO_p\,=\,\OOO'_p$ pour $p\neq 2$.
On a trivialement $\OOO_2\,=\,\OOO'_2$ si $2$ divise $d$, 
par unicit\'e de l'ordre maximal dans $A \otimes_\QQ \QQ_2$.
Nous pouvons donc supposer 
$\OOO_2 \,= \,\M_2(\ZZ_2)$ et que $\ZZ_2 1+ L_2$ est le sous-espace
des matrices de trace paire.  Mais ce sous-espace engendre
$\M_2(\ZZ_2)$ comme anneau, on en d\'eduit $\OOO_2 \subset \OOO'_2$, puis
$\OOO_2=\OOO'_2$ par maximalit\'e de $\OOO$. D'o\`u $\OOO=\OOO'$.

\medskip 
Montrons la surjectivit\'e.  Soit $L$ un \'el\'ement de $\mathcal{R}(d)$.  Montrons
d'abord que l'espace quadratique $L\otimes_\ZZ \QQ$ est isom\'etrique
\`a $A_0$. Par le th\'eor\`eme de Hasse-Minkowski, il suffit de
montrer qu'il est de d\'eterminant $2$ (dans
$\QQ^\times/(\QQ^\times)^2$), et anisotrope sur $\QQ_p$ si et
seulement si $p$ divise $d$. Le premier point est clair car $L$ est de
d\'eterminant $2d^2$.  De plus, il suffit de v\'erifier le second
point pour les premiers impairs.  En effet, par la formule du produit 
pour le symbole de Hilbert
(voir par exemple \cite[\S 4]{Serre70}) une forme quadratique sur
$\QQ^3$ suppos\'ee non d\'eg\'en\'er\'ee et d\'efinie positive est
anisotrope sur $\QQ_p$ pour un nombre impair de $p$.

Supposons donc $p$ impair. Si $p$ ne divise pas $d$, alors $L_p$ est
de d\'eterminant dans $\ZZ_p^\times$ et de rang $\geq 3$, donc
isotrope sur $\QQ_p$.  Supposons que $p$ divise $d$.  Alors par
hypoth\`ese, $\res L_p$ est anisotrope et de rang $2$ sur $\ZZ/p\ZZ$.
En particulier, le produit scalaire n'est pas identiquement nul sur $L
\otimes_\ZZ (\ZZ/p\ZZ)$.  Il existe donc $e$ dans $L_p$ avec
$\frac{1}{2} \,e \cdot e\in\ZZ_p^\times$ et on a $L_p \,=\, \ZZ_p\, e
\operp P$ avec $P$ de rang $2$ et de m\^eme r\'esidu que $L_p$.
Nous avons donc $P^\sharp = \frac{1}{p} P$.  Ainsi, il existe une
$\ZZ_p$-base de $L_p$ dans laquelle la forme quadratique $v \mapsto
\frac{v.v}{2}$ est de la forme $(x,y,z) \mapsto ax^2+\,pf(y,z)$ avec
$a\in\ZZ_p^\times$ et $f :\ZZ_p^2\ra\ZZ_p$ quadratique anisotrope
modulo $p$. Une telle forme \`a $3$ variables est manifestement
anisotrope sur $\QQ_p$.

Nous avons montr\'e que $L$ se plonge isom\'etriquement dans $A_0$.
Nous pouvons donc supposer $L \subset A_0$, le produit scalaire de
$L$ \'etant $\tr(\overline{x}y)$.  En particulier, puisque $L$ est
pair, $\n(x) \in \ZZ$ pour tout $x$ dans $L$.  Mais pour tous les $a,b$
dans $A_0$, nous avons $a^2 = - \n(a)$ puis $ab+ba = - \n(a+b)+\n(a)+\n(b)$.
Par cons\'equent, si $\epsilon_1,\epsilon_2,\epsilon_3$ est une
$\ZZ$-base de $L$, alors 
$$
\OOO_1\, =\, \ZZ\, +\, \ZZ \,\epsilon_1
\,+\, \ZZ \,\epsilon_2 \,+\, \ZZ \,\epsilon_3 \,+\, \ZZ \,\epsilon_1
\epsilon_2 \,+\, \ZZ \,\epsilon_1 \epsilon_3 \,+\, \ZZ \,\epsilon_2
\epsilon_3 \,+ \,\ZZ \,\epsilon_1 \epsilon_2 \epsilon_3
$$ 
est un sous-anneau de $A$, et donc un ordre de $A$.  Soit $\OOO$ un
ordre maximal contenant $\OOO_1$.  Alors $L(\OOO)=\OOO \cap A_0$
contient $L$, et a m\^eme d\'eterminant d'apr\`es le premier
paragraphe de la d\'emonstration. D'o\`u $L(\OOO)=L$, ce qui montre la
surjectivit\'e.  
\cqfd

\bcoro \label{cor:resLani} Soit $L\in\Rcal(d)$. Le r\'esidu de $L_2$
est isomorphe \`a $(\ZZ/2\ZZ)^3$ muni de la forme quadratique $(x,y,z)
\mapsto \frac{1}{4}(x^2+y^2+z^2) \bmod \ZZ$ si $d$ est pair, et \`a
$\ZZ/2\ZZ$ muni de la forme quadratique $x \mapsto -\frac{1}{4} x^2
\bmod \ZZ$ sinon.  En particulier, $\res L$ est anisotrope.  
\ecoro

\dem D'apr\`es la proposition pr\'ec\'edente, nous pouvons supposer
 $L=L(\OOO)$, avec $\OOO$ un ordre maximal de $A$. Notons $\OOO'$
l'ordre de Hurwitz.  Le $\ZZ_2$-r\'eseau $\OOO_2$ muni de $x\mapsto
\n(x)$ est isom\'etrique \`a $\OOO'_2$ muni de cette m\^eme forme si
$d$ est pair, et \`a $\M_2(\ZZ_2)$ muni de $x \mapsto \det(x)$ sinon.
Le sous-r\'eseau de trace nulle de $\OOO'_2$ est $\ZZ_2 \,i
\operp\ZZ_2 \,j\operp \ZZ_2\, k$ avec $\n(i)=\n(j)=\n(k)=1$.  Celui de
$\M_2(\ZZ_2)$ est $\ZZ_2 \,d \operp P$ o\`u $d$ est la matrice
diagonale $(-1,1)$, v\'erifiant $\det(d) =-1$, et $P$ est le plan
hyperbolique des matrices antidiagonales.  Le r\'esultat en d\'ecoule.
\cqfd

\medskip
Introduisons maintenant un second ensemble de r\'eseaux euclidiens
associ\'es aux ordres maximaux de $A$.  Pour tout premier $p$ impair
et $a\in\ZZ$, notons {\scriptsize $\Big(\displaystyle \frac{a}{p}
  \Big)$} le symbole de Legendre de $a$ modulo $p$.  Soit $\Scal(d)$
l'ensemble des classes d'isom\'etrie des r\'eseaux euclidiens entiers
pairs $M$ de dimension $3$ et de d\'eterminant $2d$, tels que pour
tout premier $p$ impair divisant $d$, le r\'esidu de $M_p$ soit
isomorphe au groupe $\ZZ/p\ZZ$ muni de la forme quadratique $x \mapsto
\frac{a}{p}x^2 \bmod \ZZ$, avec $a\in\ZZ$ non nul modulo $p$ tel que
{\scriptsize $\Big({\displaystyle\frac{a}{p}}\Big)=-\Big(\displaystyle
\frac{-d/p}{p}\Big)$}. Le r\'esultat suivant dit que la condition
portant sur les $\res M_p$ est superflue si $d$ est premier.

\blemm 
\label{autoscald}
Si $d$ est premier, alors tout r\'eseau euclidien entier pair $M$ de
dimension $3$ et de d\'eterminant $2d$ appartient \`a $\Scal(d)$.
\elemm

Rappelons (voir par exemple \cite[Chap.~5, \S 2 \& \S 8]{Scharlau85})
que pour tout ${\rm qe}$-{\it module} $(V,q)$, la {\it somme de Gauss}
de $(V,q)$ est 
$$ 
\ga(V,q)=|V|^{-\frac{1}{2}}\sum_{x\in V} e^{2i\pi\,q(x)}\;.  
$$ 
La somme du Gauss d'une somme orthogonale finie de ${\rm qe}$-modules 
$V_i$ est le produit des sommes de Gauss des $V_i$ (cela s'applique en particulier
\`a la d\'ecomposition en composantes primaires). La
formule de la signature de Milgram implique  $\ga(\res M) =
e^{2i\pi s/8}$ pour tout r\'eseau euclidien entier pair $M$ de rang
$s$. En particulier, cette somme de Gauss vaut $e^{3i\pi/4}$ pour $M$
de rang $3$.

\medskip

\dem Nous pouvons supposer $d$ impair.  Par la formule de Milgram et
la d\'ecomposition orthogonale $\res M \,= \,\res M_2 \operp\res M_d$,
nous avons
\begin{equation}\label{eq:calcsomgauss}
e^{3i\pi/4}\,=\,\ga(\res M) \,= \,\ga(\res M_2)\,\ga(\res M_d).
\end{equation}
Le r\'esidu de $M_2$ est isom\'etrique \`a $\ZZ/2\ZZ$ muni de la forme
quadratique $x \mapsto \epsilon \frac{x^2}{4} \bmod \ZZ$ pour un
certain signe $\epsilon=\pm 1$.  Un calcul imm\'ediat donne $\ga(\res
M_2) = e^{i\pi \epsilon /4}$.  Le r\'esidu de $M_d$ est isom\'etrique
\`a $\ZZ/d\ZZ$ muni de la forme $x\mapsto a x^2/d \bmod \ZZ$ avec
$a\in\ZZ -d\ZZ$.  D'apr\`es un th\'eor\`eme de Gauss (voir par exemple
\cite[Ch.~2]{Davenport00}), nous avons donc $\gamma(\res M_d)=$
{\scriptsize $\Big({\displaystyle \frac{a}{d}} \Big)$} pour $d \equiv
1 \bmod 4$, et $\gamma(\res M_d)=$ {\scriptsize $\Big({\displaystyle
    \frac{a}{d}}\Big)$}$i$ sinon.  Si $d\equiv 1\mod 4$, la formule
\eqref{eq:calcsomgauss} entra\^ine donc {\scriptsize
  $\Big({\displaystyle\frac{a}{d}}\Big)$} $=\epsilon=-1$. De m\^eme,
si $d\equiv 3\mod 4$ alors {\scriptsize $\Big({\displaystyle
\frac{a}{d}} \Big)$} $=\epsilon=1$. Dans les deux cas, nous avons bien
{\scriptsize $\Big({\displaystyle\frac{a}{d}}\Big)$}$=(-1)^{\frac{d+1}{2}} 
= $ {\scriptsize $ -\Big({\displaystyle\frac{-1}{d}}\Big)$}.  
\cqfd

\blemm \label{cor:resMani} 
Si $M$ est dans $\Scal(d)$ alors $\res M$ est anisotrope et $\res M_2 \simeq
\ZZ/4\ZZ$.  
\elemm

\dem Comme $\res M_p$ est anisotrope pour $p$ impair par d\'efinition,
$\res M$ est anisotrope si, et seulement si, $\res M_2$ (qui est
d'ordre $4$) l'est.  Montrons que ce dernier est isomorphe \`a
$\ZZ/4\ZZ$, ou de mani\`ere \'equivalente, que sa forme quadratique
$q$ n'est pas \`a valeurs dans $(\frac{1}{4}\ZZ)/\ZZ$.  Pour tout
premier $p$ impair divisant $d$, la somme de Gauss $\ga(\res M_p)$ est
(toujours par le th\'eor\`eme de Gauss) une racine $4$-\`eme de l'unit\'e.
Par la formule de Milgram et la multiplicativit\'e de la somme de Gauss,
$\ga(\res M_2)$ est, comme $e^{3i\pi/4}$, une racine $8$-\`eme
primitive de l'unit\'e.  En particulier, $\ga(\res M_2)$ n'est pas
r\'eelle, donc $q$ n'est pas \`a valeurs dans $(\frac{1}{2}\ZZ)/\ZZ$.
Si $q$ \'etait \`a valeurs dans $(\frac{1}{4}\ZZ)/\ZZ$, alors $\res M_2$
serait isomorphe \`a $(\ZZ/2\ZZ)^2$ muni de la forme $(x,y) \mapsto
\frac{1}{4}(u x^2+v y^2) \bmod \ZZ$, pour certains signes $u$ et $v$. Comme
la somme de Gauss de $\ZZ/2\ZZ$ muni de $x \mapsto \frac{\epsilon}{4}
x^2 \bmod \ZZ$ avec $\epsilon=\pm 1$ vaut $e^{2i\pi \epsilon /8}$, la
somme de Gauss de $\res M_2$ serait une racine $4$-\`eme de l'unit\'e :
une contradiction.  
\cqfd

\medskip
Soient $L$ un r\'eseau euclidien entier pair et $d$ un entier $\geq
1$, suppos\'es pour l'instant quelconques.  Consid\'erons les
r\'eseaux euclidiens $\Lambda=\Lambda(L;d)$ et $M = {\rm M}(L;d)$
d\'efinis par

\begin{equation}
\label{definitionM}
\Lambda = \{ x \in (d^{-1} L) \cap L^\sharp\;:\;
d\, x \cdot x \equiv 0 \bmod 2\}
\hspace{5mm} \text{et} \hspace{5mm} M \,= \,\sqrt{d} \,\Lambda\;.
\end{equation}
Autrement dit, $M$ est le plus grand sous-r\'eseau pair\footnote{Si
  $L$ est un r\'eseau entier, observer que l'application $L \ra
  \ZZ/2\ZZ$, d\'efinie par $x \mapsto x \cdot x \bmod 2$, est un
  morphisme de groupes. Son noyau est donc un r\'eseau : c'est le plus
  grand sous-r\'eseau pair de $L$. } du r\'eseau entier $N \cap
N^\sharp$ avec $N = \frac{1}{\sqrt{d}} L$. En particulier, $M$ est
entier. Par d\'efinition, nous avons $L \subset \Lambda \subset
L^\sharp$ et 
\begin{equation}\label{lambdasurl} 
\Lambda/L = \{ x \in \res L\;:\; d\,x \,=\, 0 \,\,\;\text{et}\;\,\,d\, q(x)\,=\,0\}\;.
\end{equation}

\blemm \label{lemmecurieux} Soient $L$ un r\'eseau euclidien entier pair et $d\geq 1$ un entier. 
Soit $W$ le sous-groupe des \'el\'ements $x$ de ${\rm res}\, L$ v\'erifiant $dx\,=\,0$ et $dq\,(x)\,=\,0$.
Si $q$ ne s'annule pas sur $W \cap W^\perp -\{0\}$, on a l'\'egalit\'e ${\rm M}({\rm M}(L;d);d)=L$.
\elemm

\dem 
Posons $\Lambda = \Lambda(L;d)$ et $M={\rm M}(L;d)$. Nous avons
d\'ej\`a vu $L \subset \Lambda \subset L^\sharp$ et $\Lambda/L =
W$.  Pour des raisons g\'en\'erales, nous avons alors $L \subset
\Lambda^\sharp \subset L^\sharp$, et $\Lambda^\sharp/L$ est
l'orthogonal de $W$ dans $\res L$.  Par l'hypoth\`ese sur
$W$, le seul sous-espace isotrope de $\res L$ contenu dans $W \cap
W^\perp$ est $\{0\}$.  Ainsi, $L$ est le plus grand sous-r\'eseau pair
de $\Lambda \cap \Lambda^\sharp$. Mais par d\'efinition ${\rm M}(M;d)$ est le plus
grand sous-r\'eseau pair de $N \cap N^\sharp$ o\`u $N =
\frac{1}{\sqrt{d}} M = \Lambda$ : nous avons montr\'e ${\rm
  M}(M;d)=L$.  
\cqfd

\medskip

 Notons que l'hypoth\`ese du lemme portant sur $W$ est 
automatiquement satisfaite si $\res L$ est anisotrope.
\smallskip

\bprop \label{cardeuxO}
L'application $L\mapsto M(L;d)$ induit des bijections $\Rcal(d)
\rightarrow \Scal(d)$ et $\Scal(d) \rightarrow \Rcal(d)$ qui sont
inverses l'une de l'autre.  
\eprop

\dem  
Tout \'el\'ement $L$ de $\Rcal(d)$ ou $\Scal(d)$ a son r\'esidu
anisotrope d'apr\`es le corollaire \ref{cor:resLani} et le lemme
\ref{cor:resMani}.  On en d\'eduit ${\rm M}({\rm M}(L;d);d)=L$ d'apr\`es le
lemme \ref{lemmecurieux}.  De plus, si $L$ et $L'$ sont des r\'eseaux
entiers pairs d'un espace euclidien $E$, et si $g \in {\rm O}(E)$
envoie $L$ sur $L'$, alors  $g({\rm M}(L;d))={\rm M}(L';d)$ pour
tout $d\geq 1$.  Autrement dit, l'application $L \mapsto M(L;d)$ passe
aux classes d'isom\'etrie, et il ne reste qu'\`a voir qu'elle
\'echange $\Rcal(d)$ et $\Scal(d)$.

Supposons $L$ dans $\Rcal(d)$ ou $\Scal(d)$. Posons $\Lambda=
\Lambda(L;d)$ et $M={\rm M}(L;d)$. Montrons que $\Lambda$ est d'indice
$2$ dans $L^\sharp$. Consid\'erons pour cela le sous-espace $W=
\Lambda/L$ de $\res L$, \'egalement donn\'e par la formule \eqref{lambdasurl}.
Pour tout premier $p$, notons $W_p = W \cap \res L_p$ la composante
$p$-primaire de $W$.  Soit $p$ un premier impair divisant $d$.  Alors
par les d\'efinitions de $\Rcal(d)$ et $\Scal(d)$, l'entier $p$ annule
$\res L_p$ et la forme quadratique de $\res L_p$ est \`a valeurs dans
$(\frac{1}{p}\ZZ)/\ZZ$, et donc $W_p\,=\, \res L_p$. Il ne reste
qu'\`a voir que $W_2$ est d'indice $2$ dans $\res L_2$.  Si $d$ est
impair, $\res L_2$ est isomorphe \`a $\ZZ/2\ZZ$ et sa forme
quadratique prend la valeur $\pm \frac{1}{4}$, donc $W_2=0$.
Supposons $d$ pair.  Si $L$ est dans $\Rcal(d)$, en identifiant $\res
L_2$ \`a $(\ZZ/2\ZZ)^3$ comme dans le corollaire \ref{cor:resLani},
nous constatons que $W_2$ est le sous-espace des $(x,y,z)\in
(\ZZ/2\ZZ)^3$ v\'erifiant $x+y+z = 0$. Enfin, si $L$ est dans $\Scal(d)$,
et en identifiant $\res L_2$ \`a $\ZZ/4\ZZ$ muni de la forme
quadratique $x \mapsto \frac{u}{8} x^2 \bmod \ZZ$ pour un $u \in \ZZ$
impair comme dans le lemme \ref{cor:resMani}, nous constatons que $W_2$
est le sous-espace $2 \ZZ/4\ZZ$. Dans tous ces cas, $W_2$ est bien
d'indice $2$ dans $\res L_2$.

Nous venons de montrer que $\Lambda$ est d'indice $2$ dans $L^\sharp$.
Donc $\det \Lambda = [L^\sharp:\Lambda]^2\det L^\sharp= \frac{4}{\det
  L}$. \'Ecrivons $\det L = 2 d^r$, avec $r=2$ ou $r=1$ selon que $L$
est dans $\Rcal(d)$ ou $\Scal(d)$.  Alors le d\'eterminant de $M$ est
$d^3 \det\Lambda = d^3 \frac{4}{\det L} = 2 d^{3-r}$, comme affirm\'e
dans l'\'enonc\'e.  Fixons $p$ premier impair divisant $d$ et
examinons de plus pr\`es les liens entre $\res L_p$ et $\res M_p$.
Comme $\Lambda$ est d'indice $2$ dans $L^\sharp$, nous avons
$\Lambda_p = L_p^\sharp$.  Le groupe $\res L_p$ est isomorphe \`a
$(\ZZ/p\ZZ)^r$.  Nous pouvons donc toujours trouver une
d\'ecomposition orthogonale en somme de deux $\ZZ_p$-sous-r\'eseaux
\begin{equation}\label{eq:lab} 
L_p = A_p \operp B_p
\end{equation}
avec $\det A_p \in \ZZ_p^\times$ et $B_p$ de $\ZZ_p$-rang $r$
v\'erifiant $B_p^\sharp = p^{-1} B_p$.
En particulier, $\Lambda_p = L_p^\sharp = A_p \oplus p^{-1} B_p$.
Mais nous avons une isom\'etrie\footnote{Si $(V,q)$ est un espace
  quadratique et si $m$ est dans $\ZZ$, notons $\langle m \rangle
  \otimes V$ l'espace quadratique $(V, mq)$.} $M_p \simeq \langle d
\rangle \otimes \Lambda_p$.  Comme $\langle d \rangle \otimes
B_p^\sharp$ est de d\'eterminant dans $(d/p)^r \, \ZZ_p^\times \subset \ZZ_p^\times$, car $\det
B_p^\sharp=\frac{1}{\det B_p}$ est dans $\frac{1}{p^r}\ZZ_p^\times$,
nous en d\'eduisons une isom\'etrie
\begin{equation}\label{eq:resmpap} 
\res M_p \simeq  \res (\langle d \rangle \otimes A_p) \simeq 
\langle d/p \rangle \otimes (A_p/p A_p),
\end{equation}
o\`u $A_p/p A_p$ ($\simeq (\ZZ/p\ZZ)^{3-r}$) est muni de la forme
quadratique $x \mapsto \frac{1}{p} \frac{x \cdot x}{2} \bmod \ZZ_p$.
La forme bilin\'eaire d'un ${\rm qe}$-module $W$ qui est un
$\ZZ/p\ZZ$-espace vectoriel avec $p$ premier peut \^etre vue \`a
valeurs dans $\ZZ/p\ZZ$ via l'isomorphisme naturel
$(\frac{1}{p}\ZZ)/\ZZ = \ZZ/p\ZZ$ induit par la multiplication par
$p$, et poss\`ede donc un d\'eterminant (ou ``discriminant'') qui est
un \'el\'ement de $(\ZZ/p\ZZ)^\times$ modulo les carr\'es : nous le
noterons $\delta(W)$.
La relation \eqref{eq:resmpap} entra\^ine $\delta(\res M_p) = \delta (
\langle d/p \rangle \otimes (A_p/p A_p)) = (d/p)^{3-r}
\delta(A_p/pA_p)$.  La relation \eqref{eq:lab} entra\^ine
$2p^r \equiv \det A_p \, \det B_p$ modulo les carr\'es de $\ZZ_p^\times$. 
En utilisant les congruences $\det A_p \equiv \delta( A_p/pA_p)$  et $p^{-r} \det B_p \equiv  \delta( \res B_p)$ 
dans $(\ZZ/p\ZZ)^\times$ modulo les carr\'es, on en d\'eduit $2 (d/p)^r \equiv \delta (A_p/pA_p) \, \delta( \res B_p)$.
L'iso\-mor\-phis\-me $\res L_p \simeq \res B_p$ entra\^ine donc au final
$$
2d/p\, \equiv \,\delta(\res M_p) \,\delta( \res L_p)\;.
$$ 
(toujours modulo les carr\'es de $(\ZZ/p\ZZ)^\times$.)
Supposons maintenant $L$ dans $\Rcal(d)$, et donc $r=2$.  Le
discriminant d'un plan quadratique anisotrope sur $\ZZ/p\ZZ$ \'etant
diff\'erent de $-1$, nous avons $\delta( \res L_p) \not \equiv -1$, donc
$\res M_p$ est de rang $1$ sur $\ZZ/p\ZZ$ avec $\frac{1}{2} \delta(\res
M_p) \not \equiv - d/p$. Si ${\rm res}\, M_p$ est le $\ZZ/p\ZZ$-espace vectoriel $\ZZ/p\ZZ$ 
muni de la forme quadratique $x \mapsto \frac{1}{p} a x^2 \bmod \ZZ$, le symbole de Legendre
de $\frac{1}{2} \delta(\res M_p)$ est par d\'efinition celui de $a$. Nous avons donc montr\'e
{\scriptsize $\Big({\displaystyle\frac{a}{p}}\Big)=-\Big(\displaystyle
\frac{-d/p}{p}\Big)$} : $M$ est dans $\Scal(d)$.
De m\^eme, si $L$ est dans $\Scal(d)$, on a $r=1$ et
$\frac{1}{2}\delta( \res L_p) \not \equiv -d/p$, puis $\res M_p$ est de rang
$2$ sur $\ZZ/p\ZZ$ avec $\delta(\res M_p) \not\equiv -1$, et $M$ est dans
$\Rcal(d)$.  
\cqfd

\medskip

\brema  
\label{genresd}
{\it {\rm (}Genre de $\Scal(d)${\rm )} }
Soit $M$ un r\'eseau euclidien entier pair de dimension $3$ et d\'eterminant $2d$, 
avec $d$ sans facteur carr\'e. Soit $p$ premier impair divisant $d$. 
Dans la terminologie 
de Conway \cite[Chap. 15 \S 7]{ConSlo88}, la classe d'isomorphisme
du $\ZZ_p$-r\'eseau $M_p$ est caract\'eris\'ee par son {\it symbole $p$-adique}, de la forme
$1^{2e_p} \,p^{e'_p}$ pour certains signes $e_p, e'_p\in \{\pm 1\}$. Par d\'efinition, on a une d\'ecomposition orthogonale 
$M_p = A_p \operp B_p$ avec $A_p$ de $\ZZ_p$-rang $2$ et de d\'eterminant dans $\ZZ_p^\times$, 
$B_p$ de $\ZZ_p$-rang $1$ et de d\'eterminant dans $p \ZZ_p^\times$, et 
$e_p$ (resp. $e'_p$) est le symbole de Legendre de $\det A_p$ (resp. $p^{-1}\det B_p$) modulo $p$.  
On en d\'eduit la relation
$$\Big(\displaystyle \frac{2d/p}{p}\Big) = e_p e'_p.$$ 
Par d\'efinition, le r\'eseau $M$ est dans $\Scal(d)$ si, et seulement si,
on a l'\'egalit\'e $e'_p = -\Big(\displaystyle
\frac{-2d/p}{p}\Big)$ pour tout premier $p$ impair divisant $d$ (la pr\'esence 
du $2$ dans cette formule s'explique par le passage 
de la forme quadratique \`a la forme bilin\'eaire). 
De mani\`ere \'equivalente, $M$ est dans $\Scal(d)$ si, et seulement si,
on a l'\'egalit\'e  $e_p = -\Big(\displaystyle
\frac{-1}{p}\Big)$ pour tout premier $p$ impair divisant $d$.
\erema

\section{Sur les ordres maximaux de diff\'erente principale}
\label{subsect:ordmaxdiffprinc}

Dans cette derni\`ere partie, motiv\'ee par la proposition \ref{prop:realgadeO}, 
nous donnons des caract\'erisations des 
ordres maximaux dont la diff\'erente est principale, puis 
de nombreux exemples de tels ordres.

\bprop \label{prop:tdnpstricpositf}
Toute alg\`ebre de quaternions $A$ sur $\QQ$, qui est d\'efinie, admet au
moins une classe de conjugaison d'ordres maximaux dont la diff\'erente est
principale.
\eprop

\dem Montrons tout d'abord le lemme suivant.

\blemm  \label{lem:carremoinsdiscr}
Toute alg\`ebre de quaternions $A$ sur $\QQ$, qui est d\'efinie, admet un
\'el\'ement de carr\'e $-D_A$, unique \`a conjugaison pr\`es.
\elemm

\dem
L'existence \'equivaut \`a demander qu'il existe un plongement de
$\QQ$-alg\`ebres de $\QQ(\sqrt{-D_A})$ dans $A$. 
Un tel plongement existe car les diviseurs premiers de $D_A$, et la place r\'eelle, 
sont ramifi\'es dans $\QQ(\sqrt{-D_A})$. L'unicit\'e d\'ecoule du
th\'eor\`eme de Skolem-Noether (voir \cite[p.~6]{Vigneras80}).
\cqfd

\medskip 
Soit $x\in A$ tel que $x^2=-D_A$.  Alors $\n(x)=D_A$.  De plus, $x$
\'etant entier sur $\ZZ$, il existe des ordres maximaux de $A$
contenant $x$.  Si $\OOO$ est un tel ordre, alors $\OOO x$ est un
id\'eal \`a gauche entier de $\OOO$ de norme $D_A$.  Par unicit\'e, il
est \'egal \`a la diff\'erente de $\OOO$, qui est donc principale.
\cqfd

\medskip

Illustrons la proposition \ref{prop:tdnpstricpositf} dans le cas de l'alg\`ebre de quaternions
de discriminant $D_A =11$. Il est bien connu que c'est la $\QQ$-alg\`ebre engendr\'ee $A$ par des \'el\'ements
$i$ et $j$ v\'erifiant $i^2=-1$, $j^2=-11$ et $ij=-ji$ (voir \cite[p.~98]{Vigneras80}). Elle contient exactement
deux classes de conjugaison d'ordres maximaux : voir la table de
\cite[p.~154]{Vigneras80} ou le tableau final de cette
note. V\'erifions que ces deux classes poss\`edent des repr\'esentants
$\OOO$ et $\OOO'$ contenant tous les deux l'\'el\'ement $j$, de
carr\'e $-11$, et donc sont tous les deux de diff\'erente principale. 

\medskip
\begin{itemize}
\item
D'une part, si $t=\frac{1+j}{2}$, l'ordre $\OOO=\ZZ[t]+i\,\ZZ[t]$ est de discriminant $11$,
donc maximal, et contient  $j=2t-1$. 

\medskip

\item
D'autre part, si $t'=-\frac{1}{2}+\frac{i+k}{4}$ alors $\OOO'' =\ZZ[t']+j\,\ZZ[t']$ est un ordre de $A$ (noter que $(t')^2=
-t'-1$ et $j$ normalise $\ZZ[t']$) contenant $j$. Si $\OOO'$ d\'esigne un ordre maximal  contenant $\OOO''$, alors
$\OOO'$ est non conjugu\'e \`a $\OOO$. En effet, puisque l'anneau $\ZZ[t']$ contient les $6$ unit\'es $\pm1,\pm t',\pm
(t')^2$, il n'est pas contenu dans un conjugu\'e de $\OOO$, qui ne
contient que les $4$ unit\'es $\pm1,\pm i$.
\end{itemize}
\medskip
Voici le r\'esultat de caract\'erisation des ordres maximaux de
diff\'erente principale. Deux de ces caract\'erisations feront intervenir les 
bijections des propositions \ref{carunO} et \ref{cardeuxO}.
On rappelle que si $\OOO$ est un ordre maximal de $A$, alors 
$L(\OOO)$ d\'esigne le $\ZZ$-r\'eseau euclidien des quaternions purs de
  $\OOO$ muni de $\tr(\overline{x}y)$ ; nous posons aussi $M(\OOO)={\rm M}(L(\OOO);D_A)$ (voir la formule \eqref{definitionM}).

\btheo \label{theo:carac1}
Soient $A$ une alg\`ebre de quaternions sur $\QQ$, qui est d\'efinie, et
$\OOO$ un ordre maximal de $A$. Les propri\'et\'es suviantes sont \'equivalentes~:
\begin{enumerate}
\item[(1)] la diff\'erente $\M$ de $\OOO$ est principale;
\item[(2)] l'ordre maximal $\OOO$ contient un \'el\'ement de norme (r\'eduite)
  $D_A$;
\item[(3)] l'ordre maximal $\OOO$ contient un \'el\'ement de carr\'e $-D_A$;
\item[(4)] le $\ZZ$-r\'eseau euclidien $L(\OOO)$ contient un \'el\'ement $x$ tel que $x\cdot x=2D_A$ et $x\cdot
  y=0\mod D_A$ pour tout $y\in L(\OOO)$;
\item[(5)] le $\ZZ$-r\'eseau euclidien $M(\OOO)$ contient un \'el\'ement $x$ tel que $x\cdot x=2$.
\end{enumerate}
\etheo

\dem Puisque $\M$ est l'unique id\'eal de norme $D_A$, il est principal
si, et seulement si, $\OOO$ contient un \'el\'ement de norme $D_A$ : 
les assertions (1) et (2) sont \'equivalentes.

Il est \'evident que l'assertion (3) implique l'assertion (2).
Montrons la r\'eciproque. Supposons
d'abord $D_A=2,3$.  Alors $\OOO$ est principal car $h_A=1$,
l'unicit\'e \`a conjugaison pr\`es de $\OOO$ (puisque $t_A=1$) et le
lemme \ref{lem:carremoinsdiscr} montrent que $\OOO$ contient un
\'el\'ement de carr\'e $-D_A$ (il serait bien s\^ur facile d'exhiber
un tel \'el\'ement dans ces cas).  Si $D_A\neq 2,3$, on conclut par le
lemme \ref{lem:carremoinsdiscrbis} ci-dessous.

Montrons que l'assertion (3) implique l'assertion (4). 
Soit $x \in \OOO$ v\'erifiant $x^2 = -D_A$. 
Un tel $x$ n'est pas dans $\QQ$, il est donc de trace nulle et de norme $D_A$. 
Cela montre $x \in L(\OOO)$ et $x \cdot x = 2 D_A$. De plus, $x$ est dans $\M$ par 
l'\'equivalence de (1) et (3). Nous en d\'eduisons $\tr(x \OOO) \subset D_A \ZZ$ 
par la formule \ref{dieselloc}, puis $x \cdot y \equiv 0 \bmod D_A$ pour tout $y$ dans $L(\OOO)$. 

Montrons que l'assertion (4) implique l'assertion (5). Soit $x$ dans 
$L(\OOO)$ v\'erifiant $x \cdot x = 2 D_A$ et $x \cdot y \equiv 0 \bmod D_A$
pour tout $y$ dans $L(\OOO)$. D'apr\`es la formule \ref{definitionM}, 
cela entra\^ine $\frac{1}{D_A} \, x \in \Lambda(L(\OOO),D_A)$, puis 
$x'=\frac{1}{\sqrt{D_A}} x \in \M(\OOO)$. On conclut car $x' \cdot x' = 2$.

Montrons que l'assertion (5) implique l'assertion (2). 
Soit $x$ dans $\M(\OOO)$ v\'erifiant $x' \cdot x' =2$. 
L'inclusion \'evidente $\sqrt{D_A} \M(\OOO) \subset L(\OOO)$
entra\^ine que l'\'el\'ement $x= \sqrt{D_A} x'$, de norme $D_A$, est dans
$L(\OOO)$, et donc dans $\OOO$. 
 \cqfd

\medskip

\blemm \label{lem:carremoinsdiscrbis}
Si $D_A\neq 2,3$ et si $x\in \OOO$ est de norme $D_A$, alors
$x^2=-D_A$. 
\elemm

\dem L'\'el\'ement $x$ n'appartient pas \`a $\QQ$, car $D_A$ est sans
facteur carr\'e, donc son polyn\^ome minimal est $X^2- tX +d$ avec $t
= \tr(x) \in \ZZ$ et $d=\n(x)=D_A$. Comme $A$ ramifie sur $\RR$, nous
avons $t^2 < 4d$.  De m\^eme, pour tout premier $p$ divisant $d$,
comme $A$ ramifie sur $\QQ_p$, le polyn\^ome $X^2 -tX +d$ est
irr\'eductible sur $\QQ_p$, donc $t^2-4d$ n'est pas un carr\'e dans
$\ZZ_p$. Comme $1+4 p \ZZ_p$ est constitu\'e de carr\'es de
$\ZZ_p^\times$, la factorisation $t^2-4d=t^2(1-4d/t^2)$ entra\^ine que
tout diviseur premier de $d$ divise \'egalement $t$.  Puisque $d$ est
sans facteur carr\'e, on montr\'e que $d$ divise $t$.  En utilisant
l'in\'egalit\'e $t^2 < 4d$, on en d\'eduit $t=0$, ou $t=\pm d$ et
$d<4$. \cqfd

\medskip

\medskip\noindent
{\bf Exemples. } Le plus petit discriminant d'une alg\`ebre de
quaternions sur $\QQ$, d\'efinie et ayant au moins $1$ (respectivement
$2$) classe(s) de conjugaison d'ordres maximaux dont la diff\'erente est
non principale, est $37$ (respectivement $67$).  En effet, le tableau
suivant donne, pour tous les entiers positifs $d \leq 100$ sans
facteur carr\'e ayant un nombre impair de facteurs premiers,

$\bullet$~ le nombre $t(d)$ de classes de conjugaison d'ordres maximaux
dans une alg\`ebre de quaternions sur $\QQ$ d\'efinie et de discriminant
$d$ (voir \cite[p.~152]{Vigneras80} pour une formule exacte), ainsi que

$\bullet$~ le nombre $t_{dnp}(d)$ (strictement inf\'erieur \`a $t(d)$ par la
proposition \ref{prop:tdnpstricpositf}) de classes de conjugaison
d'ordres maximaux dont la diff\'erente est non principale.

D'apr\`es les propositions \ref{carunO} et \ref{cardeuxO}, 
$t(d)$ est aussi le nombre de classes d'\'equivalence de 
formes quadratiques ternaires enti\`eres  
d\'efinies positives de d\'eterminant $2d$ 
appartenant au genre d\'ecrit dans la remarque \ref{genresd}.
Nous utilisons les tables de formes ternaires de Brandt et Intrau, 
recalcul\'ees et rendues disponibles sur le site de Nebe et Sloane \cite{NebSlo} par Schiemann.
Dans la terminologie de ces tables, 
le {\it discriminant} d'une telle forme d\'esigne l'entier $-d$. 
Nous en d\'eduisons par inspection la ligne $t(d)$ de la table ci-dessous.

De plus, l'\'equivalence entre les assertions (1) et (5) du th\'eor\`eme \ref{theo:carac1}
montre que $t_{dnp}(d)$ est le nombre de classes d'\'equivalence de formes ternaires ci-dessus
qui ne repr\'esentent pas l'entier $1$.  \'Etant donn\'e que dans les tables sus-cit\'ees
les formes ternaires sont donn\'ees sous forme r\'eduite, une telle forme repr\'esente $1$ si, et seulement si,
son premier coefficient est $1$  (alternativement, on peut aussi v\'erifier en utilisant par exemple le logiciel SAGE (algorithme LLL)
que le r\'eseau euclidien associ\'e  \`a cette forme a ses plus courts vecteurs de carr\'e scalaire \'egal \`a $2$).
On en d\'eduit la ligne $t_{dnp}(d)$ de la table ci-dessous.

$$
\begin{array}{|c||c|c|c|c|c|c|c|c|c|c|c|c|c|c|c|}
\hline d & 2 & 3 & 5 & 7 & 11 & 13 & 17 & 19 & 23 & 29 &  
30  
& 31 & 37 & 41& 42  
  \\\hline
t     & 1 & 1 & 1 & 1 & 2 & 1 & 2 & 2 & 3 & 3 & 1 & 3 & 2 & 4  & 1 
\\\hline
t_{dnp}& 0 & 0 & 0 & 0 & 0 & 0 & 0 & 0 & 0 & 0 & 0 & 0 & 1 & 0  & 0  
\\\hline
\hline
d & 43 & 47 & 53  & 59 & 61  &66 &   67 & 70 & 71 & 73 & 
 78  & 79 & 83 & 
89
& 97 \\\hline
t    & 3 & 5 & 4  & 6 & 4 &  2  & 4 & 1 & 7 & 4 & 1 & 6 & 7 & 7 & 5 \\\hline
t_{dnp} & 1 & 0 & 1 & 0 & 1 &  0 & 2 & 0 & 0 & 2 & 0 & 1 & 1 & 1 & 3 \\\hline 
\end{array}
$$

\bigskip 

Dans l'article \cite{Ibukiyama82}, Ibukiyama donne une formule 
pour le nombre de classes de conjugaison
d'ordre maximaux de $A$ contenant un \'el\'ement de carr\'e $-d$ avec $d=D_A$, 
ou ce qui revient au m\^eme, 
pour la quantit\'e $t(d)-t_{dnp}(d)$ d'apr\`es le th\'eor\`eme \ref{theo:carac1}. 
Dans le cas o\`u $d$ est un nombre premier impair $p$, 
cette formule est particuli\`erement simple et due \`a Deuring. Elle s'\'ecrit $$t(p)-t_{dnp}(p) = \frac{h(-p) + h(-4p)}{2},$$
o\`u $h(-m)$ d\'esigne le nombre de {\it classes d'\'equivalence propre 
de formes quadratiques binaires enti\`eres positives et primitives de discriminant $-m$}  (voir la remarque 2.13 dans  \cite{Ibukiyama82}). 
Cette formule confirme la table ci-dessus.
Nous pourrions en fait la red\'emontrer sans grande difficult\'e 
\`a partir de l'\'equivalence entre les assertions (1) et (5) du th\'eor\`eme \ref{theo:carac1},
et des propositions \ref{carunO} et \ref{cardeuxO} 
(observer, en guise de point  de d\'epart, 
que pour $M$ dans $\Scal(p)$ et $x$ dans $M$ avec $x \cdot x = 2$, 
l'orthogonal de $\ZZ x$ dans $M$ est de dimension $2$ et de d\'eterminant \'egal \`a $p$ ou $4p$).

{\small \bibliography{biblio} }

\bigskip
{\small
\noindent 
\begin{tabular}{l}
Laboratoire de math\'ematique d'Orsay,\\
UMR 8628 Univ. Paris-Sud et CNRS,\\
Universit\'e Paris-Saclay,\\
91405 ORSAY Cedex, FRANCE\\
{\it e-mail: gaetan.chenevier@math.cnrs.fr}
\end{tabular}
\hfill 
\begin{tabular}{l}
Laboratoire de math\'ematique d'Orsay,\\
UMR 8628 Univ. Paris-Sud et CNRS,\\
Universit\'e Paris-Saclay,\\
91405 ORSAY Cedex, FRANCE\\
{\it e-mail: frederic.paulin@math.u-psud.fr}
\end{tabular}
}

\end{document}

%% file: macros_fr_pdf.tex
\usepackage[T1]{fontenc}
\usepackage[utf8]{inputenc}
\usepackage{amsmath,amssymb,epsfig,bbm}
\usepackage{stmaryrd,stix}
\usepackage[colorlinks]{hyperref}
\usepackage{comment}
\usepackage{color}
\usepackage{textcomp}

\usepackage[acronym,toc]{glossaries}
\usepackage{glossary-longragged}
\usepackage{glossary-superragged}
\newglossary[slg]{symbols}{sym}{sbl}{List of Symbols}

\usepackage[textsize=small]{todonotes}
\usepackage{varwidth}

\usepackage{url}

\usepackage{enumitem}
\setlist[itemize]{labelindent=*, leftmargin=.5 truecm,nosep}

\usepackage[nottoc]{tocbibind}

\usepackage{fancyhdr}
\usepackage[us,12hr]{datetime} 
\fancypagestyle{plain}{
\fancyhf{}
\rfoot{ {\ddmmyyyydate\today}}
\lfoot{\thepage}
}

\usepackage{amsthm}

\theoremstyle{plain}
\newtheorem{defi}{D\'efinition}
\newtheorem{prop}[defi]{Proposition}
\newtheorem{theo}[defi]{Th\'eor\`eme}
\newtheorem{conj}[defi]{Conjecture}
\newtheorem{lemm}[defi]{Lemme}
\newtheorem{coro}[defi]{Corollaire}

\theoremstyle{definition}
\newtheorem{rema}[defi]{Remarque}
\newtheorem{exem}[defi]{Exemple}
\newtheorem{exems}[defi]{Exemples}

\newcommand{\bdefi}{\begin{defi}}
\newcommand{\edefi}{\end{defi}}
\newcommand{\bprop}{\begin{prop}}
\newcommand{\eprop}{\end{prop}}
\newcommand{\btheo}{\begin{theo}}
\newcommand{\etheo}{\end{theo}}
\newcommand{\btheofr}{\begin{theofr}}
\newcommand{\etheofr}{\end{theofr}}
\newcommand{\blemm}{\begin{lemm}}
\newcommand{\brema}{\begin{rema}}
\newcommand{\erema}{\end{rema}}
\newcommand{\bexer}{\begin{exem}}
\newcommand{\eexer}{\end{exem}}
\newcommand{\bexem}{\begin{exem}}
\newcommand{\eexem}{\end{exem}}
\newcommand{\bexems}{\begin{exems}}
\newcommand{\eexems}{\end{exems}}
\newcommand{\bconj}{\begin{conj}}
\newcommand{\econj}{\end{conj}}
\newcommand{\elemm}{\end{lemm}}
\newcommand{\bcoro}{\begin{coro}}
\newcommand{\ecoro}{\end{coro}}
\newcommand{\dem}{\noindent{\bf Démonstration. }}

\usepackage{mathrsfs}
\renewcommand\mathcal{\mathscr}

\newcommand{\E}{{\cal E}}

\newcommand{\G}{{\cal G}}

\newcommand{\M}{{\cal M}}
\newcommand{\N}{{\cal N}}
\newcommand{\OOO}{{\cal O}}

\newcommand{\Q}{{\cal Q}}

\newcommand\Rcal{\mathcal R}
\newcommand{\Scal}{{\mathcal S}}

\newcommand{\X}{{\mathcal X}}
\newcommand{\Y}{{\mathcal Y}}

\newcommand{\maths}[1]{{\mathbb #1}}  

\newcommand{\FF}{\maths{F}}
\newcommand{\HH}{\maths{H}}

\newcommand{\NN}{\maths{N}}

\newcommand{\QQ}{\maths{Q}}
\newcommand{\RR}{\maths{R}}

\newcommand{\ZZ}{\maths{Z}}

\newcommand{\ga}{\gamma}

\newcommand{\ra}{\rightarrow}

\newcommand{\covol}{\operatorname{Covol}}

\newcommand{\cqfd}{\hfill$\Box$}

\newcommand{\gengeod}
{\operatorname{\widecheck{\G\,}\!\!}}

\newcommand{\n}{\operatorname{\tt n}}

\newcommand{\res}{\operatorname{res}}

\newcommand{\tr}{\operatorname{\tt tr}}

\newcommand{\udh}{\operatorname{U}_2(\HH)}
\newcommand{\Vol}{\operatorname{Vol}}

\newcommand{\GL}{\operatorname{GL}}

\newcommand{\SL}{\operatorname{SL}}

\usepackage{interval}

\usepackage{mathtools}
\makeatletter
\DeclareRobustCommand\widecheck[1]{{\mathpalette\@widecheck{#1}}}
\def\@widecheck#1#2{%
    \setbox\z@\hbox{\m@th$#1#2$}%
    \setbox\tw@\hbox{\m@th$#1%
       \widehat{%
          \vrule\@width\z@\@height\ht\z@
          \vrule\@height\z@\@width\wd\z@}$}%
    \dp\tw@-\ht\z@
    \@tempdima\ht\z@ \advance\@tempdima2\ht\tw@ \divide\@tempdima\thr@@
    \setbox\tw@\hbox{%
       \raise\@tempdima\hbox{\scalebox{1}[-1]{\lower\@tempdima\box
\tw@}}}%
    {\ooalign{\box\tw@ \cr \box\z@}}}
\makeatother

\newcounter{fig}



%% file: fig_racines16.pdf_t
\begin{picture}(0,0)%
\includegraphics{fig_racines16.pdf}%
\end{picture}%
\setlength{\unitlength}{3812sp}%
\begingroup\makeatletter\ifx\SetFigFont\undefined%
\gdef\SetFigFont#1#2#3#4#5{%
  \reset@font\fontsize{#1}{#2pt}%
  \fontfamily{#3}\fontseries{#4}\fontshape{#5}%
  \selectfont}%
\fi\endgroup%
\begin{picture}(1516,1514)(2843,-2978)
\put(3421,-2401){\makebox(0,0)[lb]{\smash{{\SetFigFont{11}{13.2}{\rmdefault}{\mddefault}{\updefault}{\color[rgb]{0,0,0}$x\mapsto x^{-1}$}%
}}}}
\put(3151,-1996){\makebox(0,0)[lb]{\smash{{\SetFigFont{11}{13.2}{\rmdefault}{\mddefault}{\updefault}{\color[rgb]{0,0,0}d'ordre $4$}%
}}}}
\end{picture}%